\documentclass[12pt,a4paper,english]{article}

\usepackage{amsthm}
\usepackage{epsfig}
\usepackage{alltt}
\usepackage{makeidx}
\usepackage{newlfont}
\usepackage{amsmath}
\usepackage{amssymb}
\usepackage{amsfonts}
\usepackage{amscd}

\makeindex

\newcommand{\op}{\operatorname}
\newcommand{\korr}{ }

\newtheorem{theorem}{Theorem}[section]

\newtheorem{Theorem}[theorem]{Theorem}

\newtheorem{Proposition}[theorem]{Proposition}

\newtheorem{Corollar}[theorem]{Corollary}
\newtheorem{lemma}[theorem]{Lemma}
\newtheorem{Lemma}[theorem]{Lemma}
\newcommand{\BLemma}{\begin{lemma}}
\newcommand{\ELemma}{\end{lemma}}
\newtheorem{Satz}[theorem]{Theorem}

\theoremstyle{definition}
\newtheorem{definition}[theorem]{Definition}
\newcommand{\BDefinition}{\begin{definition}}
\newcommand{\EDefinition}{\end{definition}}

\newtheorem{Definition}[theorem]{Definition}
\theoremstyle{remark}

\newtheorem{Behauptung}[theorem]{Behauptung}

\newcommand{\BBehauptung}{\begin{Behauptung}}
\newcommand{\EBehauptung}{\end{Behauptung}}
\newtheorem{Bemerkung}[theorem]{Bemerkung} 
 
\newtheorem{Remark}[theorem]{}
\newtheorem{Remarkc}[theorem]{Remark}
\newcommand{\BBemerkung}{\begin{Bemerkung}}
\newcommand{\EBemerkung}{\end{Bemerkung}}
\newcommand{\BProposition}{\begin{Proposition}}
\newcommand{\EProposition}{\end{Proposition}}
\newcommand{\BSatz}{\begin{Satz}}
\newcommand{\ESatz}{\end{Satz}}
\newtheorem{Bemerkungen}[theorem]{Bemerkungen}
\newcommand{\BBemerkungen}{\begin{Bemerkungen}}
\newcommand{\EBemerkungen}{\end{Bemerkungen}}
\newcommand{\BTheorem}{\begin{Theorem}}
\newcommand{\ETheorem}{\end{Theorem}}
\newtheorem{Beispiel}[theorem]{Beispiel}

\newcommand{\BBeispiel}{\begin{Beispiel}}
\newcommand{\EBeispiel}{\end{Beispiel}}
\newtheorem{Beispiele}[theorem]{Beispiele}
\newcommand{\BBeispiele}{\begin{Beispiele}}
\newcommand{\EBeispiele}{\end{Beispiele}}

\newcommand{\DC}{{\Bbb C}}

\newcommand{\DZ}{{\Bbb Z}}

\newcommand{\ra}{\rightarrow}
\newcommand{\lra}{\longrightarrow}
\newcommand{\sra}{\twoheadrightarrow}

\newcommand{\da}{\downarrow}

\newcommand{\hra}{\hookrightarrow}

\newcommand{\IFF}{\Leftrightarrow}

\newcommand{\sira}{\stackrel{\sim}{\rightarrow}}

\input{xy}
\xyoption{all}
\begin{document}
\title{Andersen Filtration and Hard Lefschetz}
\author{Wolfgang Soergel}

%\address{Universit\"at Freiburg\\Mathematisches Institut
%\\Eckerstra\ss e 1\\ D-79104 Freiburg\\Germany}
%\email{soergel@@mathematik.uni-freiburg.de}
%\thanks{partially supported by the 
%TMR Algebraic Lie Theory ERB FMRX-CT97-0100}
\maketitle 
\begin{center}
  \emph{For Joseph Bernstein}
  \end{center}

\begin{abstract}
On the space of homomorphisms from a Verma module to an indecomposable tilting
module of the BGG-category $\cal{O}$
we define a natural filtration following Andersen \cite{AFil}
and establish a formula expressing the dimensions of the filtration
steps in terms of coefficients of Kazhdan-Lusztig polynomials. 
\end{abstract}
%%%%%%%%%%%%%%%%%%%%%%%%%%%%%%%%%%%%%%%%%%%%%%%%%%%%%%%%%%%%%%%%%%%%%%%%%%%%%%%%%%%%%

\section{Introduction}
Indecomposable tilting modules in category $\mathcal{O}$ were classified by
Collingwood and Irving \cite{CI} well before this terminology existed under the
name of selfdual Verma flag modules. More precisely, they proved 
that applying the 
indecomposable projective functors of Bernstein-Gelfand to simple 
instead of projective
Verma modules, you get precisely the 
indecomposable selfdual Verma flag modules instead of
the indecomposable projectives, and that these indecomposable 
selfdual Verma flag modules nowadays called tilting modules are
classified by their highest weight.
Now we can define a filtration on the space of homomorphisms 
from a Verma module to a tilting
module by analogy of what Andersen \cite{AFil} did in the 
algebraic group case. The main result
of this article is a description of the dimensions of 
the subquotients of this filtration in terms
of Kazhdan-Lusztig polynomials.

To be more precise, let me introduce some notation.
Let $\mathfrak{g} \supset \mathfrak{b} \supset \mathfrak{h}$ be a 
semisimple complex Lie algebra, a Borel
and a Cartan. Let $\rho\in \mathfrak{h}^\ast$ 
be the halfsum of roots from $\mathfrak{b}$
and let $\DC[\mathbb{C} \rho] =T$ denote
the regular functions on the line $\mathbb{C} \rho \subset\mathfrak{h}^\ast.$ 
This is a quotient of
$S\frak{h}^\ast = \DC[\mathfrak{h}] ,$ and
every linear form $v:\mathbb{C} \rho\ra \mathbb{C} $ defines 
an isomorphism  with a polynomial ring
$\mathbb{C} [v] \overset{\sim}{\rightarrow}T.$ 
For a weight $\lambda \in \mathfrak{h}^\ast$ we 
form the Verma module $\Delta (\lambda) 
=U (\mathfrak{g})\otimes_{U(\mathfrak{b})} \mathbb{C}_\lambda 
\in \mathfrak{g}\op{-mod}$ and the deformed Verma
module
\begin{equation*}
\Delta_T (\lambda) = U(\mathfrak{g}) \otimes_{U(\mathfrak{b})} 
(\mathbb{C}_\lambda \otimes T) \in \mathfrak{g}\op{-mod-}T
\end{equation*}
Here and henceforth tensor products without any specification are to be
understood over 
  $\DC$.
The $T$-action is meant to only move the last 
tensor factor, however the $\mathfrak{b}$-action on $\mathbb{C}_\lambda
\otimes T$ comes via the obvious surjection 
$\mathfrak{b} \twoheadrightarrow \mathfrak{h}$ from
the $\mathfrak{h}$-action given by the tensor product 
action $H (a \otimes f) = \lambda (H) a \otimes f
+a \otimes Hf$ for $a \in \mathbb{C}_\lambda = \mathbb{C}$ and $f \in T$.
Starting with 
the deformed Verma and taking
the $T$-dual ``weight space by weight space'' and 
twisting the $\mathfrak{g}$-action on the result
with a Chevalley automorphism we also get a deformed dual 
Verma module $\nabla_T (\lambda) \in \mathfrak{g}\op{-mod-}T.$
The universal properties of Verma modules will then lead to a 
canonical embedding
\begin{equation*}
\op{can} : \Delta_T (\lambda) \hookrightarrow \nabla_T (\lambda)
\end{equation*}
which gives an isomorphism between the (analogues of the) 
highest weight spaces and is in fact
a basis of the $T$-module $\op{Hom}_{\mathfrak{g}-T} 
(\Delta_T (\lambda), \nabla_T (\lambda))$.
The Jantzen filtration can be understood as the
filtration of our Verma $\Delta (\lambda)$ by the images 
of the $\op{can}^{-1} 
(\nabla_T (\lambda) v^i) $ for $i = 0,1,2, \ldots$
under the natural projection $\Delta_T (\lambda) 
\twoheadrightarrow \Delta (\lambda).$
Next let $\nu \in \mathfrak{h}^\ast$ be such that 
the Verma module $\Delta (\nu)$ is simple and let
$E \in \mathfrak{g}\op{-mod}$ be finite dimensional. 
Then $E \otimes \Delta (\nu)$ is tilting
 and we may consider
the composition pairing
\begin{equation*}
\op{Hom} (\Delta_T (\lambda), E\otimes 
\Delta_T (\nu)) \times \op{Hom}
(E \otimes \Delta_T (\nu), \nabla_T (\lambda)) \rightarrow 
\op{Hom} (\Delta_T (\lambda), 
\nabla_T (\lambda))
\end{equation*}
where homomorphisms are understood in the category of
$\mathfrak{g}$-$T$-bimodules.
As we remarked already, the pairing essentially lands in $T$. 
Furthermore we will prove that
the paired spaces actually are free of finite rank over $T$, 
thus our pairing can be rewritten as a map,
actually an embedding
\begin{equation*}
\op{Hom} (\Delta_T (\lambda), E\otimes 
\Delta_T (\nu)) \hookrightarrow \op{Hom}
(E \otimes \Delta_T (\nu), \nabla_T (\lambda))^\ast
\end{equation*}
with the $\ast$ meaning a $T$-dual.
Andersen's filtration is defined by taking on the right side of this embedding
the filtration obtained by multiplying with the $v^i$ from the right,
then taking the preimage of this filtration under our embedding, 
and finally the image of this preimage 
under the projection onto $\op{Hom}_{\mathfrak{g}} (\Delta (\lambda), 
E \otimes \Delta (\nu))$
specializing $v$ to $0$ alias applying $\otimes_T \mathbb{C}$.

In this paper we explain how to calculate the dimensions of 
the subquotients $\bar{F}^{i} $ of this Andersen filtration 
on $\op{Hom}_{\mathfrak{g}} (\Delta (\lambda), E \otimes \Delta (\nu))$.
More precisely, we identify the dimensions of the subquotients 
of the induced filtration on $\op{Hom}_{\mathfrak{g}}
(\Delta (\lambda), K)$ for $K \subset E \otimes \Delta (\nu)$ an 
indecomposable direct summand with
coefficients of Kazhdan-Lusztig polynomials $P_{y,x}(q)$  
as they are introduced in \cite{KL-C}.
And to be completely explicit,  the general formula we prove 
as Theorem \ref{MT} means 
 in the principal block for arbitrary $x,y$ in the Weyl group the formula
\begin{equation*}
\sum_{i \geq 0} \op{dim}_{\mathbb{C}} \bar{F}^{i}
\op{Hom}_{\mathfrak{g}} (\Delta (-y\rho-\rho),
K)\; q^{(l(x)-l(y)-i)/2}= P_{y,x}(q)
\end{equation*}
in case $K$ has highest weight $(-x\rho-\rho)$ and thus is the 
indecomposable tilting module $K =K (-x\rho-\rho)$
with this highest weight. 

The  proof given in the last section
proceeds roughly speaking by translation to the 
Koszul-dual geometric side, where we run into the hard Lefschetz.
More precisely, the embedding  giving rise to Andersen's
filtration is identified  with the embedding of a costalk 
of the equivariant intersection cohomology complex of
a Schubert variety 
into its stalk at the same point, both understood in the  equivariant 
derived category of a point.
This identification in turn passes through identifying both sides 
with the same construction in bimodules over polynomial rings, i.e.
passing through a ``coherent picture''. More precisely, in  sections 
2-5 we explain the translation from category $\cal{O}$ to the
coherent picture, culminating in \ref{DTM}.
The  translation from geometry to the coherent picture 
is discussed thereafter.

The arguments given even show that the Andersen filtration 
coincides with the filtration on our spaces of homomorphisms
coming from the $\Bbb{Z}$-graded
structure introduced in \cite{BGSo}, although we do not make this explicit.
This statement is very similar to the 
semisimplicity of the subquotients of
the Jantzen filtration proved in \cite{BB-J}, 
but the method to obtain it is quite different.
I would like to know how to directly relate both results,
as this would give an alternative proof of the mentioned semisimplicity.

\section{Deformation of category $\cal{O}$}
\begin{Remark}
  In this and the next section we repeat results of \cite{GJ}
 in a language adapted to our goals,
which is also very close to
the language introduced in \cite{FieC}.  
 Let $S = S \frak{h}
= \DC[\frak{h}^*]$ be the symmetric algebra of $\frak{h}$.
We consider the category $\op{Kring}^S$ of all commutative unitary rings
$T$ with a distinguished morphism
$\varphi : S \ra T.$ 
%%%%%%%%%%%%%%%%%%%%%%%%%%%%%%%%%%%%%%%%%%%%%%%%%%%%%%%%%%%%%%%%%%%%%%%%%%%%%%
Given $T\in \op{Kring}^S$
  we consider the category $\frak{g}\op{-Mod_{\DC}-} T$ of all
  $\frak{g}\op{-}T$-bimodules on which the right and left actions of
  $\Bbb{C}$ coincide.  
\end{Remark}
\begin{Definition}
For $T=(T,\varphi)\in\op{Kring}^S$
we define 
in any bimodule $M \in
  \frak{g}\op{-Mod_{\DC}-}T$ for any $\lambda \in \frak{h}^{\ast}$ 
the {\bf deformed weight space}
  $M^{\lambda}$  by the formula
  $$M^{\lambda} =M^{\lambda}_T=\{ m \in M 
\mid (H - \lambda (H)) m = m \varphi (H)
  \quad\forall H \in \frak{h}\}$$
\end{Definition}
%%%%%%%%%%%%%%%%%%%%%%%%%%%%%%%%%%%%%%%%%%%%%%%%%%%%%%%%%%%%%%%%%%%%%%%%%
\begin{Remark}
For $M \in \frak{g}\op{-Mod_{\DC}-}T $ the canonical map from the
direct sum of its deformed weight spaces to
$M$ is always an injection
$\bigoplus_\lambda M^\lambda\hra M.$
For $T = \DC[\frak{h}^{\ast}]$ this is evident, since
the weight spaces $M^{\lambda}$
considered as $T \otimes T$-modules have support in the graphs of
$(\lambda +): 
\frak{h}^{\ast} \ra \frak{h}^{\ast}$ and these graphs are pairwise disjoint.
In general our weight spaces have support in the preimage of our graphs
under the map 
$\op{Spec} (\DC[\frak{h}^{\ast}] \otimes
T) \ra \op{Spec} (\DC[\frak{h}^{\ast}]
\otimes\DC[\frak{h}^{\ast}])$
induced by
$\op{id} \otimes \varphi$
and thus are disjoint as well.
\end{Remark}
%%%%%%%%%%%%%%%%%%%%%%%%%%%%%%%%%%%%%%%%%%%%%%%%%%%%%%%%%%%%%%%%%%%%%%%%%%%%%%%
\begin{Definition}
For every $T\in\op{Kring}^S$  we define in our
category of bimodules a full subcategory, the
{\bf deformed category}
  $$\cal{O}(T) \subset \frak{g}\op{-Mod_{\DC}-}T$$
as the category of all bimodules $M$ which are locally finite for
$\frak{n}=[\frak{b},\frak{b}]$
and decompose as the direct sum $M=\bigoplus M^\lambda$ 
of their deformed weight spaces.
\end{Definition}
\begin{Remark}
Prominent objects of this category are the  {\bf deformed
Verma modules}
  $$\Delta_{T} (\lambda) =
  \op{prod}_{\frak{b}}^\frak{g}(\Bbb{C}_{\lambda}\otimes T) =U (\frak{g})
  \otimes_{U(\frak{b})} (\Bbb{C}_{\lambda}\otimes T)$$
for $\lambda \in
  \frak{h}^{\ast}$, where it is understood that the right action of
$T$ acts  only on the last tensor factor, whereas the left action 
comes from the left action of $U(\frak{b})$ on
  $\Bbb{C}_{\lambda} \otimes T$ which we get via the canonical surjection
  $\frak{b} \twoheadrightarrow \frak{h}$ from the tensor
action of  $\frak{h},$
  where $H \in \frak{h}$
acts  on $\Bbb{C}_{\lambda}$ 
via the scalar  $\lambda (H)$ and on $T$ by multiplication with 
$\varphi (H)$.
\end{Remark}
\begin{Remark}
The category $\cal{O} (T)$ is stable under tensoring from
the left with finite dimensional representations of $\frak{g},$ where
as left action of $\frak{g}$ on such a tensor product we understand the tensor
action and as right action of
$T$ its right action on the second tensor factor. Along with
a bimodule $\cal{O} (T)$ also contains all its subquotients.  In case
$T = \Bbb{C}$ and $\varphi$  the evaluation at the zero of $
\frak{h}^{\ast},$ the category $\cal{O} (T)$ 
specializes up to some missing finiteness
conditions to the usual category
$\cal{O}$ of Bernstein-Gelfand-Gelfand, and
$\Delta_\DC(\lambda)=\Delta(\lambda)$ is the
Verma module with highest weight
$\lambda.$
\end{Remark}
%%%%%%%%%%%%%%%%%%%%%%%%%%%%%%%%%%%%%%%%%%%%%%%%%%%%%%%%%%%%%%%%%%%%%%%%%%%%%%
\begin{Definition}
We now consider the opposed Borel of $\frak{b}$ with respect to $\frak{h}$
to be denoted $\bar{\frak{b}} \subset \frak{g}$ and
for $\lambda\in \frak{h}^\ast$ consider the subbimodule
$$\nabla_{T} (\lambda) \subset
  \op{ind}_{\bar{\frak{b}}}^\frak{g}(\Bbb{C}_{\lambda}\otimes T)
  =\op{Hom}_{U(\bar{\frak{b}})} (U (\frak{g}), \Bbb{C}_{\lambda}\otimes T)$$
defined as the sum of all deformed weight spaces of the
$\op{Hom}$-space in question.
We call it the {\bf deformed Nabla-module}
of highest weight
 $\lambda.$
\end{Definition}
\begin{Remark}
Under the identification given by restriction of our
$\op{Hom}$-spaces with $\op{Hom}_{\Bbb{C}}(U (\frak{n}), \Bbb{C}_{\lambda}
\otimes T)$ our $\nabla_{T}(\lambda)$ corresponds to those homomorphisms,
which are different from zero on at most
finitely many $\frak{h}$-weight spaces of $U(\frak{n})$.
The deformed nablas also belong to
$\cal{O} (T).$
\end{Remark}
%%%%%%%%%%%%%%%%%%%%%%%%%%%%%%%%%%%%%%%%%%%%%%%%%%%%%%%%%%%%%%%%%%%%%%%%%%%%%%%%
\begin{Remark}
All weight spaces of $\nabla_{T} (\lambda)$
and $\Delta_{T} (\lambda)$ are free over
$T$ and finitely generated, and if $T$ is not zero, the deformed
weight spaces of weight $(\lambda - \nu)$
in both modules have the rank
$\op{dim}_{\Bbb{C}}U(\frak{n} )^{\nu}.$ 
%%%%%%%%%%%%%%%%%%%%%%5
We have canonical morphisms
$T \overset{\sim}{\ra} \Delta_{T}(\lambda)^{\lambda}
\hookrightarrow \Delta_{T} (\lambda)$
and $\nabla_{T}(\lambda)
  \twoheadrightarrow \nabla_{T}(\lambda)^{\lambda} \overset{\sim}{\ra} T$ of
$T$-modules and for any extension $T \ra T^{\prime}$
canonical isomorphisms
$\Delta_{T} (\lambda) \otimes_{T} T^{\prime} \overset{\sim}{\ra}
  \Delta_{T^{\prime}} (\lambda)$ and  $\nabla_{T} (\lambda) \otimes_{T}
  T^{\prime} \overset{\sim}{\ra} \nabla_{T^{\prime}} (\lambda).$
\end{Remark}
%%%%%%%%%%%%%%%%%%%%%%%%%%%%%%%%%%%%%%%%%%%%%%%%%%%%%%%%%%%%%%%%%%%%%%%%%%%%%%
\begin{Remark}
We now choose for our Lie algebra an involutive automorphism 
$\tau: \frak{g} \ra \frak{g}$
with $\tau |_{\frak{h}} = -\op{id}$ and define a contravariant functor
$$d= d_{\tau} : \frak{g}\op{-Mod_{\DC}-}T \ra \frak{g}\op{-Mod_{\DC}-}T$$
by letting $dM \subset \op{Hom}_{-T} (M,T)^{\tau}$ be the sum of all
deformed weight spaces in the space of homomorphisms of right 
$T$-modules from $M$ to $T$  with its contragredient
$\frak{g}$-action twisted
by $\tau$.  If $M\in
\frak{g}\op{-Mod_{\DC}-}T $ is the sum of its deformed weight spaces,
we have a canonical morphism
$M\ra ddM,$ and if in addition all deformed weight spaces of
$M$ are free and finitely generated over
$T,$ this canonical  morphism is an isomorphism.
\end{Remark}
%%%%%%%%%%%%%%%%%%%%%%%%%%%%%%%%%%%%%%%%%%%%%%%%%%%%%%%%%%%%%%%%%%%%%%%
\begin{Remark}
The restriction onto the highest deformed
weight space defines together with the universal
property of the induced representation a canonical homomorphism
$$
  \op{Hom}_{-T} (\op{prod}^{\frak{g}}_{\bar{\frak{b}}} (\Bbb{C}_{-\lambda}
  \otimes T),T) \ra \op{ind}_{\bar{\frak{b}}}^{\frak{g}} \op{Hom}_{-T}
  (\Bbb{C}_{\lambda} \otimes T,T)
$$
and considering the deformed weight spaces we see that
it induces an isomorphism of bimodules 
$$d\Delta_{T} (\lambda) \sira \nabla_{T}
(\lambda).
$$
With our preceding remarks we also get $d\nabla_{T}
(\lambda) \cong \Delta_{T} (\lambda) .$ By the tensor identity,
i.e. since tensoring with a representation
of a Lie algebra commutes with tensor-inducing a
representation from a subalgebra, furthermore
$E \otimes \Delta_{T} (\lambda)$ admits a filtration with subquotients
$\Delta_{T} (\lambda + \nu),$ where $\nu$ runs over the multiset $ P (E)$ of
weights of $E.$ Since $E \otimes ?$ commutes up to the choice
of an isomorphism $dE\cong  E$
with our duality
$d$,
we deduce an analogous result for
$E\otimes  \nabla_{T} (\lambda).$
\end{Remark}
%%%%%%%%%%%%%%%%%%%%%%%%%%%%%%%%%%%%%%%%%%%%%%%%%%%%%%%%%%%%%%%%%%%%%%%%%%%%
\begin{Proposition}\label{EHo}
  \begin{enumerate}
  \item For all $\lambda$ the restriction to the deformed weight space of 
$\lambda$ together with the two canonical identifications
$\Delta_{T}(\lambda)^{\lambda} \overset{\sim}{\ra} T$ and
$\nabla_{T} (\lambda)^{\lambda} \overset{\sim}{\ra}T$ induces an isomorphism
$$\op{Hom}_{\cal{O}(T)} (\Delta_{T} (\lambda),
      \nabla_{T}(\lambda)) \overset{\sim}{\ra} T.$$
\item For $\lambda \neq \mu$ in $\frak{h}^{\ast}$
we have
$\op{Hom}_{\cal{O} (T)} (\Delta_{T} (\lambda),
\nabla_{T} (\mu))=0.$
\item For all $ \lambda, \mu \in\frak{h}^{\ast}$ we have
$\op{Ext}^{1}_{\cal{O} (T)} (\Delta_{T} (\lambda), \nabla_{T} (\mu))=0 .$ 
\end{enumerate}
\end{Proposition}
\begin{proof} We prove (3), the simpler case of 
spaces of homomorphisms is treated in the same way.
Let $R^+\subset \frak{h}^\ast$ denote the roots of $\frak{n} $
and $|R^+\rangle\subset \frak{h}^\ast$ the submonoid generated by $R^+$
and $\leq$ the partial order on
$\frak{h}^\ast$ with $\lambda\leq\mu\IFF \mu \in \lambda + |R^+\rangle.$
Every short exact sequence $\nabla_{T} (\mu)
\hookrightarrow M \twoheadrightarrow \Delta_{T} (\lambda)$
with $M \in \cal{O} (T)$ and $ \lambda \not\leq \mu$
splits, since any preimage in
$M^{\lambda}$ of the canonical generator of $\Delta_{T} (\lambda)$
already is annihilated by $\frak{n} $ and thus induces a splitting.
In case $\lambda \leq \mu$ we use our duality $d$ to pass to the dual
situation. This proves the triviality of the extension in question.
\end{proof}
%%%%%%%%%%%%%%%%%%%%%%%%%%%%%%%%%%%%%%%%%%%%%%%%%%%%%%%%%%%%%%%%%%%%%%%%%%%%%%%%
\begin{Corollar}\label{BCH}
Let  $M,N \in \cal{O} (T).$
If $M$ is a direct summand of an object with finite $\Delta_T$-flag and
$N$ a direct summand of an object with finite $\nabla_T$-flag, then
the space of homomorphism $\op{Hom}_{\cal{O}(T)} (M,N)$ is a
finitely generated projective $T$-module and for any
ring extension $T \ra T^{\prime}$ the obvious
map defines an isomorphism
$$\op{Hom}_{\cal{O}(T)} (M,N) \otimes_{T} T^{\prime} 
\;\;\overset{\sim}{\ra}\;\; \op{Hom}_{\cal{O}(T^{\prime})}
(M \otimes_{T} T^{\prime}, N\otimes_{T} T^{\prime})$$
\end{Corollar}
\begin{proof}
This follows directly from \ref{EHo} by induction on the lengths 
of the flags. 
\end{proof}
%%%%%%%%%%%%%%%%%%%%%%%%%%%%%%%%%%%%%%%%%%%%%%%%%%%%%%%%%%%%%%%%%%%%%%%%%%%%%%%%%%%%
\begin{Remark}\label{he}
If $Q\in\op{Kring}^S$ is a field and if for all roots $\alpha$ the
coroots $\alpha^{\vee}$ are not mapped to $\Bbb{Z} \subset Q$ under
$S \ra Q$, then the category
$\cal{O}(Q)$ is semisimple, (i.e.\ all surjections split)
and its simple objects are the
$\Delta_{Q} (\lambda) = \nabla_{Q} (\lambda)$ for
$\lambda \in \frak{h}^{\ast}.$
\end{Remark}
%%%%%%%%%%%%%%%%%%%%%%%%%%%%%%%%%%%%%%%%%%%%%%%%%%%%%%%%%%%%%%%%%%%%%%%%%%%%%%%%%%%%
%%%%%%%%%%%%%%%%%%%%%%%%%%%%%%%%%%%%%%%%%%%%%%%%%%%%%%%%%%%%%%%%%%%%%%%%%%%%%%%%%%%%
\section{Deforming indecomposable tilting modules}
\begin{Remark}
Let $D=S_{(0)}$ be the local ring at zero of $\frak{h}^{\ast}.$
For $\lambda \in \frak{h}^{\ast}$ with $\Delta (\lambda)$ simple
the canonical map defines an isomorphism
$$\Delta_{D} (\lambda) \overset{\sim}{\ra} \nabla_{D} (\lambda)$$
Indeed, we only need to show that this map gives isomorphisms on all
deformed weight spaces, and these are free of finite rank over the
local ring $D.$ By Nakayama's Lemma we thus only need to show
that our map becomes an isomorphism under 
$?\otimes_{D} \Bbb{C}$, and this follows
directly from the simplicity assumption on $\Delta (\lambda)$.
\end{Remark}
%%%%%%%%%%%%%%%%%%%%%%%%%%%%%%%%%%%%%%%%%%%%%%%%%%%%%%%%%%%%%%%%%%%%%%%%%%%%%%%%%%%%
\begin{Definition} 
Given $T\in\op{Kring}^S$
let $\cal{K} (T) \subset \cal{O} (T)$ denote the smallest
subcategory, which 
\begin{enumerate}
\item  contains all $\Delta_{T} (\lambda)$ for which the canonical map gives an
  isomorphism $\Delta_{T} (\lambda) \overset{\sim}{\ra} \nabla_{T} (\lambda)$,
\item   is stable under tensoring with finite dimensional representations of
  $\frak{g}$, \item  is stable under forming direct summands.
\end{enumerate}
We call $\cal{K} (T)$ the category of {\bf $T$-deformed tilting modules}.
\end{Definition}
%%%%%%%%%%%%%%%%%%%%%%%%%%%%%%%%%%%%%%%%%%%%%%%%%%%%%%%%%%%%%%%%%%%%%%%%%%%%%%%%%%%%%%
\begin{Remark}
For $\DC=\DC_0\in \op{Kring}^S$ the objects
of $\cal{K} (\DC)$
are the  tilting modules of the usual
BGG-category $\cal{O}.$
\end{Remark}
\begin{Proposition}\label{CTD}
If $T\in\op{Kring}^S$
is a complete local ring ``under $S$''
such that the preimage in $S$ of its maximal ideal is
just the vanishing ideal of the origin in
$\frak{h}^{\ast}$, then the specialization
$$?\otimes_{T} \Bbb{C} : \cal{K} (T)\ra\cal{K} (\Bbb{C})$$
induces a bijection on isomorphism classes, and under this
bijection indecomposables correspond to  indecomposables.
\end{Proposition}
%%%%%%%%%%%%%%%%%%%%%%%%%%%%%%%%%%%%%%%%%%%%%%%%%%%%%%%%%%%%%%%%%%%%%%%%%%%%%%%%%%%%%
\begin{proof}
The tilting modules from $\cal{O}$ are precisely the direct summands of
tensor products of simple Vermas with finite dimensional
representations.
\korr{All such tensor products $K\in \cal{K} (\Bbb{C})$ lift by definition.
If $K_T\in \cal{K} (T)$ is such a lift, we deduce 
from \ref{BCH} that the canonical map leads to  an isomorphism
$\DC\otimes_T\op{End} K_T\sira \op{End} K_\DC$
of finite dimensional $\DC$-algebras. General results
\cite{Ben} or 
\cite{CuRe}, I, (6.7) concerning the lifting of idempotents now show
that any projection of $K_\DC$ to a direct summand 
can be lifted to a projection of $K_T$ to a direct summand,
which gives surjectivity on isomorphism classes in our Proposition.
The same argument, now applied to an arbitrary $K\in \cal{K} (\Bbb{C}),$
shows that only indecomposable objects from $\cal{K} (T)$ can
go to indecomposable objects from $\cal{K} (\DC).$ 
Similarily, any lifting of an isomorphism has to be an isomorphism,
since every lift of a unit in an endomorphism ring has to be a unit,
and this establishes the claimed bijection on isomorphism classes.}
\end{proof}
%%%%%%%%%%%%%%%%%%%%%%%%%%%%%%%%%%%%%%%%%%%%%%%%%%%%%
\begin{Remark}
For $\lambda\in \frak{h}^\ast$ we let
$K_T(\lambda)\in \cal{K} (T)$ denote the $T$-deformation
of the indecomposable tilting module $K(\lambda)\in \cal{O}$ with
highest weight $\lambda.$
\end{Remark}
%%%%%%%%%%%%%%%%%%%%%%%%%%%%%%%%%%%%%%%%%%%%%%%%%%%%%%%%%%%%%%%%%%%%%%%%%%%%%%%%%%%%%%%%
%%%%%%%%%%%%%%%%%%%%%%%%%%%%%%%%%%%%%%%%%%%%%%%%%%%%%%%%%%%%%%%%%%%%%%%%%%%%%%%%%%%%%%%%%%
\section{The Andersen filtration}

\begin{Remark}
Fix $K \in \frak{g}\op{-Mod_{\DC}-}T$
and
$\lambda \in \frak{h}^{\ast}$.
To increase readability we use the abbreviations
$\Delta_{T} (\lambda)=\Delta,$ $\nabla_{T} (\lambda)=\nabla$
and $\op{Hom}_{\frak{g}-T}=\op{Hom}$ and
consider the
$T$-bilinear pairing
$$\op{Hom} (\Delta, K) \times \op{Hom}
(K,\nabla) \ra \op{Hom} (\Delta,
\nabla) = T$$
given by composition.
If for any $T$-module $H$ we denote by $H^\ast$
the $T$-module $\op{Hom}_T(M,T),$ then our pairing induces a map
$$E=E_\lambda(K):
\op{Hom} (\Delta , K) \ra \op{Hom}
(K, \nabla)^{\ast}$$
If $K$ is tilting, then by
\ref{BCH} our map $E$ is a map 
between finitely generated projective $T$-modules.
If in addition $T\in\op{Kring}^S$ is an integral domain and
$Q = \op{Quot}T$ satisfies the assumptions of remark \ref{he}, thus
$\cal{O}(Q)$ is semisimple with simple objects $\Delta_{Q} (\lambda) =
\nabla_{Q} (\lambda),$ then our pairing is nondegenerate over $Q$ and
our map $E_\lambda(K)$
induces an isomorphism over $Q$
and in particular is an injection.
If now $T =\Bbb{C} [[v]]$ is the ring of formal power series around 
the origin on a
line $\DC \delta \subset \frak{h}^{\ast},$ which isn't 
contained in any reflection hyperplane
of the
Weyl group, then $Q=\op{Quot}\Bbb{C} [[v]]$
satisfies our assumptions of remark
\ref{he}. 
If now $K\in \cal{K}(\Bbb{C} [[v]])$
is a deformed tilting module, we can use the embedding
$$E_\lambda(K):\op{Hom} (\Delta , K) \hra
  \op{Hom} (K, \nabla)^{\ast}$$
of free
$\Bbb{C} [[v]]$-modules of finite rank to restrict the
obvious filtration of the right hand side by the
$v^{i}\op{Hom} (K, \nabla)^{\ast}$ and thus get a filtration on
$\op{Hom} (\Delta, K)
=\op{Hom}_{\frak{g}-\Bbb{C} [[v]]} (\Delta_{\Bbb{C} [[v]]} (\lambda), K).$
\end{Remark}
%%%%%%%%%%%%%%%%%%%%%%%%%%%%%%%%%%%%%%%%%%%%%%%%%%%%%%%%%%%%%%%%%%%%%%%%%%%%%%
\begin{Definition}\label{DAF}
Given $K_\DC\in\cal{K}(\DC)$ a tilting module of
$\cal{O}$ and $K\in \cal{K}(\DC[[v]])$ a $\DC[[v]]$-deformation
of $K_\DC$ in the sense of \ref{CTD} with $S\ra\DC[[v]]$
the restriction to a formal neighbourhood of the origin in the line  
$\DC\rho$ with $\rho$ as in
\ref{rhod},
the image of the filtration defined above under
specialization
$?\otimes_{\Bbb{C}[[v]]} \Bbb{C}$
will be called the {\bf Andersen-filtration on}
$\op{Hom}_{\frak{g}} (\Delta (\lambda), K_{\Bbb{C}}).$
\end{Definition}
\begin{Remark}
We leave it to the reader to show that this filtration is
independent of the choice of the deformation,
which is only well defined up to isomorphism.
The goal of this work is to determine the dimensions of the
subquotients of the Andersen filtration on
$\op{Hom}_{\frak{g}} (\Delta (\lambda), K(\mu))$
for all $\lambda,\mu\in\frak{h}^\ast$ or more precisely their description
as coefficients of Kazhdan-Lusztig polynomials.
\end{Remark}
\begin{Theorem}\label{MT}
The dimensions of the subquotients
of the  Andersen filtrations satisfy the identities
\begin{displaymath}
\sum_i\dim_\DC \bar{F}^i\op{Hom}_{\frak{g}}(\Delta (\lambda_{\bar{y}}), 
K(\lambda_{\bar{x}}))q^{(l(x)-l(y)-i)/2}= P_{y,x}(q)
\end{displaymath}
\end{Theorem}

\begin{Remark}
The proof will be given only in the last section, 
but let me explain here what all this notation means.
  We start with a $\rho$-dominant weight
  $\lambda\in\frak{h}^\ast_{\op{dom}}$ in the sense of \ref{rhod}. It gives
  two subgroups $W_{\bar{\lambda}}\supset W_\lambda$ of the Weyl group as
  explained in \ref{WLL}, and $\bar{x},\bar{y}$ denote cosets of
  $W_{\bar{\lambda}}/ W_\lambda$ with $x,y$ their longest representatives.
  Finally $\lambda_{\bar{x}}=w_{\bar{\lambda}}\bar{x}\cdot\lambda$ is to be
  understood as in \ref{NZ} with $w_{\bar{\lambda}}$ the longest element of
  $W_{\bar{\lambda}},$ and  $P_{y,x}$ is the
 Kazhdan-Lusztig polynomial with respect
  to the Coxeter group $W_{\bar{\lambda}}$ and its length function $l.$ 
  In fact the arguments given
  in this article show that the Andersen filtration coincides with the grading
  filtration induced from the graded version of $\cal{O}$, but I felt
  incapable to explain this in the framework of this article.
\end{Remark}

%%%%%%%%%%%%%%%%%%%%%%%%%%%%%%%%%%%%%%%%%%%%%%%%%%%%%%%%%%%%%%%%%%%%%%%%%%%%%%
\begin{Remark}
The Jantzen filtration on a Verma module 
$\Delta (\lambda)$ certainly induces
a filtration on
$\op{Hom}_{\frak{g}} (P(\mu), \Delta (\lambda))$ for
$P (\mu) \twoheadrightarrow \Delta (\mu)$
the indecomposable projective cover of $\Delta (\mu)$ in $\cal{O}.$
This filtration in turn comes in the same way from the embedding
$\Delta_{\Bbb{C} [[v]]}(\lambda)\ra\nabla_{\Bbb{C} [[v]]}(\lambda)$
or more precisely the embeddings
$$\op{Hom}_{\frak{g}-\Bbb{C} [[v]]}(P_{\Bbb{C} [[v]]} (\mu),
\Delta_{\Bbb{C} [[v]]} (\lambda)) \hookrightarrow
\op{Hom}_{\frak{g}-\Bbb{C} [[v]]} (P_{\Bbb{C} [[v]]}(\mu), 
\nabla_{\Bbb{C} [[v]]} (\lambda))$$
induced by them, where
$P_{\Bbb{C} [[v]]}(\mu)\sra \Delta_{\Bbb{C} [[v]]} (\mu)$
are the indecomposable projective covers in
$\cal{O}(\Bbb{C} [[v]]).$ This shows the analogy of both
filtrations. 
\korr{In fact, the contravariant equivalence  explained in
\cite{So-CT} from the category of 
Verma flag modules to itself, mapping projectives to tilting
modules, induces a map on homomorphism spaces, and this map  should 
identify both filtrations. However I cannot prove this without using the
Jantzen conjecture.}
\end{Remark}
%%%%%%%%%%%%%%%%%%%%%%%%%%%%%%%%%%%%%%%%%%%%%%%%%%%%%%%%%%%%%%%%%%%%%%%%%%%%%%%%%%%%
%%%%%%%%%%%%%%%%%%%%%%%%%%%%%%%%%%%%%%%%%%%%%%%%%%%%%%%%%%%%%%%%%%%%%%%%%%%%%%%%%%%%%%
\section{Deformed translation}
\begin{Remark}
Let $Z \subset U(\frak{g})$ be the center, so that $Z
\otimes T$ acts on any bimodule $M \in \frak{g}\op{-Mod_{\DC}-}T $.
We now consider the push-out diagram of $\Bbb{C}$-algebras
  $$\begin{array}{ccccc}
    & & Z \otimes T & &\\
    &\nearrow & & \searrow &\\
    Z \otimes \Bbb{C} [\frak{h}^{\ast}] & & & & \Bbb{C}[\frak{h}^{\ast}]
    \otimes T\\
    & \searrow & \ & \nearrow\\
    & &\Bbb{C}[\frak{h}^{\ast}] \otimes \Bbb{C} [\frak{h}^{\ast}]&&
 \end{array}$$
 %%%%%%%%%%%%%%%%%%%%%%%%%%%%%%%%%%%%%%%%%%%%%%%%%%%%%%%%%%%%%%%%%%%%%%%%%%%
where for $\xi : Z \ra \Bbb{C} [\frak{h}^{\ast}]$ we always take the variant
of the Harish-Chandra homo\-mor\-phism with $\xi (z) -z \in U \frak{n} $.
It leads to a finite ring extension and the same holds thus also for
both downward arrows of our diagram.
The graph of the addition of $\lambda \in \frak{h}^{\ast}$ is an
irreducible closed subset of
$\op{Spec} (\Bbb{C}[\frak{h}^{\ast}] \otimes
 \Bbb{C} [\frak{h}^{\ast}])$ and the same holds for its image in
$\op{Spec} (Z \otimes \Bbb{C} [\frak{h}^{\ast}]).$ 
The preimage in $\op{Spec} (Z\otimes T)$ of this image will
be denoted
$\Xi_{\lambda}\subset \op{Spec} (Z \otimes T)$.
By definition $\Delta_{T}
 (\lambda)$ and $\nabla_{T} (\lambda)$ both have 
 support in $\Xi_{\lambda}$
as
 $Z\otimes T$-modules.
\end{Remark}
%%%%%%%%%%%%%%%%%%%%%%%%%%%%%%%%%%%%%%%%%%%%%%%%%%%%%%%%%%%%%%%%%%%%%%%%%%
\begin{Lemma}
The support in
$\op{Spec} (Z \otimes T)$
of any element of an object
$M \in \cal{O} (T)$ is contained in a finite union of sets
of the form $\Xi_{\lambda}$
with $ \lambda \in \frak{h}^{\ast}.$
\end{Lemma}
\begin{proof}
Let $v$ be our element.
We may assume $v \in M^{\lambda}$
for some $\lambda \in \frak{h}^{\ast}.$
We may further assume the submodule generated by
$v$ to be contained in
$\bigoplus_{\mu \leq \nu}M^{\lambda + \mu}$
for any integral dominant weight $\nu \in X^{+}.$
The object
$$U (\frak{g}) \otimes_{U(\frak{b})} \tau_{\leq \lambda +\nu} \left(U(\frak{b})
\otimes_{U(\frak{h})} (\Bbb{C}_{\lambda} \otimes T)\right)$$
with hopefully selfexplaining $\tau_{\leq \lambda + \nu}$
has a finite $\Delta_{T}$-flag
and our $v$ is contained in the image of a homomorphism of said object
to $M$.
\end{proof}
%%%%%%%%%%%%%%%%%%%%%%%%%%%%%%%%%%%%%%%%%%%%%%%%%%%%%%%%%%%%%%%%%%%%%%%%%
\begin{Definition}\label{rhod}
Let $\rho=\rho(R^+)$ be the halfsum of positive roots.
We put
$$\frak{h}^\ast_{\op{dom}} = \{ \lambda \in \frak{h}^\ast \mid
\langle \lambda + \rho, \alpha^{\vee} \rangle \not\in \{ -1,-2,
\ldots\} \; \forall \alpha \in R^{+}\}$$
and call the elements of this set 
$\rho$-dominant weights.
We use the usual notation $w\cdot\lambda=w(\lambda+\rho)-\rho$
for the action of the Weyl group translated to the
fixed point $-\rho.$
\end{Definition}
%%%%%%%%%%%%%%%%%%%%%%%%%%%%%%%%%%%%%%%%%%%%%%%%%%%%%%%%%%%%%%%%%%%%%%%%%%%
\begin{Satz}[Decomposition of deformed categories]
Let $T$ be an $S_{(0)}$-ring, i.e.\ the
morphism $S\ra T$ should factor through the local ring $S_{(0)}$
of $\frak{h}^\ast$ at the origin.
Then we have a decomposition
$$\cal{O} (T) = \prod_{\lambda \in \frak{h}^{\ast}_{\op{dom}} }
\cal{O}_{\lambda}(T)$$
where $\cal{O}_{\lambda} (T)$ consists of all $M \in \cal{O}(T)$ which
satisfy
$M = \bigoplus_{\nu \in \lambda + \Bbb{Z} R} M^{\nu} $ and $
\op{supp}_{Z \otimes T} M
\subset \bigcup_{w \in W} \Xi_{w \cdot \lambda}.$
\end{Satz}
%%%%%%%%%%%%%%%%%%%%%%%%%%%%%%%%%%%%%%%%%%%%%%%%%%%%%%%%%%%%%%%%%%%%%%%%%%%%%
\begin{proof}
From $\Xi_{\lambda} \cap \Xi_{\mu} \neq \emptyset$ we get for
$T=S_{(0)}$ already $W \cdot
\lambda = W \cdot \mu.$
The rest of the argument can be copied from the case
$T = \Bbb{C},$ see \cite{BG}.
\end{proof}
\begin{Remark}
As in the non-deformed case we have for $\lambda,
\mu \in \frak{h}_{\op{dom}}^{\ast}$
with integral difference $\lambda -\mu \in X$ translation functors
$$T^{\mu}_{\lambda} : \cal{O}_{\lambda} (T) \ra \cal{O}_{\mu} (T)$$
which are exact, satisfy adjunctions 
$(T^{\mu}_{\lambda}, T^{\lambda}_{\mu})$
and have all the usual properties.
We call them
{\bf deformed translations}.
The category of deformed tilting modules in one of our blocks will
be denoted $\cal{K}(T)\cap \cal{O}_\lambda(T)=\cal{K}_\lambda(T).$
\end{Remark}
\begin{Remark}\label{WLL}
Let us put $D=S_{(0)}.$
If $T$ is a  $D$-algebra, then
for $\lambda \in \frak{h}^{\ast}_{\op{dom}}$
the deformed Verma module $\Delta_{T} (\lambda)$
is projective in $\cal{O}_{\lambda} (T)$. The isotropy group of a weight
$\lambda\in\frak{h}^{\ast}$  under
the dot-action of the Weyl group will be denoted $W_\lambda,$
the isotropy group of its coset $\bar{\lambda}=\lambda+\langle R\rangle$
under the root lattice will be denoted $W_{\bar{\lambda}}.$
The longest element of $W_{\bar{\lambda}}$ will be denoted
$w_{\bar{\lambda}},$
the rings of invariants for the natural actions of
$W_\lambda\subset W_{\bar{\lambda}}$ on $D$ will be denoted
$D^\lambda\supset D^{\bar{\lambda}}.$
\end{Remark}
\begin{Satz}[Deformation of projectives]
The functor $?\otimes_{D} \Bbb{C} :
\cal{O} (D) \ra \cal{O} (\Bbb{C})$ induces a bijection
between isomorphism classes of finitely generated projective
objects in both categories.
\end{Satz}
%%%%%%%%%%%%%%%%%%%%%%%%%%%%%%%%%%%%%%%%%%%%%%%%%%%%%%%%%%%%%%%%%%%%%%%%%%%%%%%
\begin{proof}
\cite{So-A}.
\end{proof}
\begin{Definition}
Given $\lambda \in \frak{h}^{\ast}$ let
$P_{D} (\lambda) \in \cal{O} (D)$ denote the finitely
generated projective specalizing to
$P (\lambda)$ under $?\otimes_{D} \Bbb{C} $. We call it the {\bf deformation of the
projective $P (\lambda)$.} Given $\lambda \in \frak{h}^{\ast}_{\op{dom}}$ 
%%%%%%%%%%%%%%%%%%%%%%%%%%%%%%%%%%%%%%%%%%%%%%%%%%%%%%%%%%%%%%%%%%%%%%%%%%%%%%%
we use for the deformed antidominant projective the abbreviation
$P_{D} (w_{\bar{\lambda}} \cdot\lambda)=A_D(\lambda)=A(\lambda).$
\end{Definition}
\begin{Satz}[Endomorphisms of antidominant projectives]
Given $\lambda \in \frak{h}^{\ast}_{\op{dom}}$ the
multiplication defines a surjection
$Z \otimes D \twoheadrightarrow 
\op{End}_{\cal{O}(D)} A (\lambda).$
If $(+\lambda)^{\sharp} : S \ra S$ denotes
the comorphism of $(+\lambda) : \frak{h}^{\ast}\ra
\frak{h}^{\ast}$, then the composition
$$Z \otimes D \overset{\xi \otimes \op{id}}{\lra} 
S \otimes D
\overset{(+\lambda)^{\sharp}\otimes\op{id}}{\lra} S 
\otimes D \ra D \otimes_{D^{\bar{\lambda}}}
D$$
has image $D^{{\lambda}} \otimes_{D^{\bar{\lambda}}} D$ and the same
kernel as the surjection considered before and we thus get an isomorphism
$$D^{{\lambda}} \otimes_{D^{\bar{\lambda}}} D \overset{\sim}{\ra}
\op{End}_{\cal{O}(D)} A ( \lambda)$$
\end{Satz}
\begin{proof}
For $\lambda$ integral the proof is given in \cite{HCH}.
The proof in general is essentially the same.
\end{proof}
\begin{Remark}
For better transparency we use hereafter frequently the notation
$\op{Hom}_{\cal{O}(D)}=\op{Hom}$ and $\op{End}_{\cal{O}(D)}=\op{End}.$
Any choice of a deformed antidominant projective
$A ( \lambda)$ for
$\lambda \in \frak{h}^{\ast}_{\op{dom}}$ gives us via the rule
$\Bbb{V} = \Bbb{V}_{D}
= \op{Hom}_{\cal{O}(D)}(A( \lambda),\;)$   
an exact functor
$$
\Bbb{V}:\cal{O} (D) \ra D^{{\lambda}}\op{-Mod_{\DC}-}D$$
which is different from zero only on $\cal{O}_{\lambda} (D)$.
If further $\mu \in \frak{h}^{\ast}_{\op{dom}}$ is given with
$\lambda -\mu \in X$ and $W_{\mu} \supset W_{\lambda}$ and if
we choose an isomorphism $T^{\lambda}_{\mu} A( \mu) \overset{\sim}{\ra}
A(\lambda)$, we get a commutative diagram
$$\begin{array}{ccc}
D^{{\mu}} \otimes_{D^{\bar{\mu}}} D& 
\overset{\sim}{\ra} & \op{End}
A(\mu)\\
\downarrow & & \downarrow \\
D^{{\lambda}} \otimes_{D^{\bar{\mu}}} D & 
\overset{\sim}{\ra} & \op{End}
A(\lambda)
\end{array}$$
with the left vertical induced from the embedding $D^{{\mu}} \subset D^{{\lambda}}$
and the right vertical given by
$T^{\lambda}_{\mu}$ and our isomorphism, see
\cite{HCH}.  If we fix such an isomorphism and in addition choose an
adjunction $(T^{\lambda}_{\mu}, T^{\mu}_{\lambda})$,
then we get isomorphisms
$$\op{Hom} (A(\mu), T^{\mu}_{\lambda} M)
  \overset{\sim}{\ra} 
\op{Hom} (T^{\lambda}_{\mu} A (
  \mu), M) \overset{\sim}{\ra} 
\op{Hom} (A( \lambda),M)$$
which lead to an isomorphism of functors, up to which the diagram
  $$\begin{array}{ccc}
    \cal{O}_{\lambda} (D) & \overset{\Bbb{V}}{\ra} & D^{{\lambda}} 
\op{-Mod_{\DC}-}D\\[2mm]
    T_{\lambda}^{\mu} \downarrow & & \downarrow \op{res}\\[2mm]
    \cal{O}_{\mu} (D) & \overset{\Bbb{V}}{\ra} & D^{{\mu}}\op{-Mod_{\DC}-}D
\end{array}$$
commutes.  Using the adjunctions we also find an isomorphism of functors,
up to which the diagram
$$\begin{array}{ccl}
  \cal{O}_{\mu} (D) 
& \overset{\Bbb{V}}{\ra} & D^{{\mu}} \op{-Mod_{\DC}-}D\\[2mm]
  T^{\lambda}_{\mu} \downarrow & 
& \;\;\;\;\;\;\;\;\downarrow D^{{\lambda}}
  \otimes_{D^{{\mu}}}?\\[2mm]
  \cal{O}_{\lambda}(D) &\overset{\Bbb{V}}{\ra} & D^{{\lambda}}
\op{-Mod_{\DC}-}
  D
\end{array}$$
commutes.
\end{Remark}
%%%%%%%%%%%%%%%%%%%%%%%%%%%%%%%%%%%%%%%%%%%%%%%%%%%%%%%%%%%%
\begin{Remark}
Given $\lambda, \mu \in \frak{h}^*_{\op{dom}}$
with integral difference one may more generally consider the
translations
$
T^{\mu}_{\lambda} : \mathcal{O}_{\lambda} (D) \ra \mathcal{O}_{\mu} (D)$
which can be written as $T^{\mu}_{\lambda}
\cong T_{\nu}^{\mu} T^{\nu}_{\lambda}$ for one and any
$\nu \in \frak{h}^*_{\op{dom}}$
with integral difference to $\lambda$ and $\mu$ and the
property $W_{\nu} = W_{\lambda}
\cap W_{\mu}$. \korr{This can be deduced from the known effects on Verma 
modules using the classification of projective functors \cite{BG}.}
If we form $D^{\mu}_{\lambda} = D^{W_{\lambda} \cap W_{\mu}} \in D^{\mu}
\op{-Mod_{\Bbb{C}}-}D^{\lambda}$, then we may interpret both diagrams as
one diagram commuting up to natural isomorphism, namely
the diagram
\begin{displaymath}
\xymatrix{
\mathcal{O}_{\lambda}(D) \ar[r]^-{\Bbb{V}} \ar[d]_{T^{\mu}_{\lambda}} 
& D^{\lambda}
\op{-Mod_{\Bbb{C}}-}D \ar[d]^{D^{\mu}_{\lambda} \otimes_{D^{\lambda}}?}\\
\mathcal{O}_{\mu} (D) \ar[r]^-{\Bbb{V}} & D^{\mu}\op{-Mod_{\Bbb{C}}-}D
}
\end{displaymath}
If we pass to the adjoints of the vertical functors, we get another diagram
commuting up to natural isomorphism, namely
\begin{displaymath}
\xymatrix{
\mathcal{O}_{\mu} (D) \ar[r]^-{\Bbb{V}} \ar[d]_{T^{\lambda}_{\mu}} &
D^{\mu}\op{-Mod_{\Bbb{C}}-}D\ar[d]^{\op{Hom}_{D^{\mu}} (D_{\lambda}^{\mu},
?) }\\
\mathcal{O}_{\lambda} (D) \ar[r]^-{\Bbb{V}} & D^{\lambda}\op{-Mod_{\Bbb{C}}-}
D
}
\end{displaymath}

\end{Remark}
%%%%%%%%%%%%%%%%%%%%%%%%%%%%%%%%%%%%%%%%%%%%%%%%%%%%%%%%%%%%%%%%%%%%%%%%%%%%%%%%%%
\begin{Theorem}[Struktursatz for deformed tilting modules]\label{DTM}
The functors $\Bbb{V}$ are fully faithful on deformed tilting
modules.
More precisely we have for any $\lambda\in\frak{h}^\ast_{\op{dom}}$:
\begin{enumerate}
\item
Given $K \in \cal{K}_\lambda (D) $
and $F \in \cal{O}_{\lambda} (D)$
an object with $\Delta_D$-flag the functor $\Bbb{V}$ induces an isomorphism
$$\op{Hom}_{\frak{g}-D} (F, K) \overset{\sim}{\ra} 
\op{Hom}_{D^{\lambda}-D} (\Bbb{V} F, \Bbb{V}K)$$
\item
Given $K \in \cal{K}_\lambda (D) $
and $F \in \cal{O}_{\lambda} (D)$ an object with a $\nabla_D$-flag
the functor
$\Bbb{V}$ induces an isomorphism
$$\op{Hom}_{\frak{g}-D} (K,F) \overset{\sim}{\ra} \op{Hom}_{D^{\lambda}-D}
(\Bbb{V} K, \Bbb{V} F)$$
\end{enumerate}
\end{Theorem}
%%%%%%%%%%%%%%%%%%%%%%%%%%%%%%%%%%%%%%%%%%%%%%%%%%%%%%%%%%%%%%%%%%%%%%
\begin{Remark}
In greater generality the first statement is proven as
Theorem 10 in \cite{FieC}: The functors  $\Bbb{V}$ are even
fully faithful on arbitrary objects with a finite
$\Delta_D$-flag.
\end{Remark}
\begin{Remark}
In the non-deformed case $T = \Bbb{C}$ the functor $\Bbb{V}$
is fully faithful on the category of tilting modules of a given block.
Indeed for any maximal ideal
$\chi \subset Z$ and arbitrary projective functors
$F, G : U / \chi U \op{-mod} \ra U \op{-mod}$ and an arbitrary
Verma module $\Delta$ with $\chi \Delta =0$
applying our functors to $\Delta$ defines a  bijection
\begin{displaymath}
\op{Trans}_{U/\chi U \op{-mod}\ra} (F,G) 
\overset{\sim}{\rightarrow} \op{Hom}_{\frak{g}}
(F\Delta, G\Delta)
\end{displaymath}
\korr{where on the left hand side we mean the set of transformations 
from the functor $F$ to the functor $G$ and did only specify the
start category, since this is the most subtle point in this business.}
%%%%%%%%%%%%%%%%%%%%%%%%%%%%%%%%%%%%%%%%%%%%%%%%%%%%%%%%%%%%%%%%%%%
For projective Vermas this is shown in
\cite{BG}, and since by \cite{BGG-Aff}
the enveloping algebra surjects onto the
$\op{ad}$-finite endomorphisms of every Verma, the proof given
there works more generally for every Verma.
The embedding of a simple Verma $\Delta_e$  into a projective Verma
$\Delta_p$ thus gives bijections
\begin{displaymath}
\op{Hom}_{\frak{g}} (F \Delta_e, G \Delta_e) 
\overset{\sim}{\rightarrow} \op{Hom}_{\frak{g}}
(F\Delta_p, G\Delta_p)
\end{displaymath}
%%%%%%%%%%%%%%%%%%%%%%%%%%%%%%%%%%%%%%%%%%%%%%%%%%%%%%%%%%%%%%%%%%%%%%%%%%%%%%
Since it also gives bijections
$\Bbb{V} F \Delta_e \overset{\sim}{\rightarrow} \Bbb{V}
F \Delta_p,$ the claim follows.
In the nondeformed case however the faithfulness on morphisms
from tiltings to dual Vermas or from Vermas to tiltings does not hold.
\end{Remark}
%%%%%%%%%%%%%%%%%%%%%%%%%%%%%%%%%%%%%%%%%%%%%%%%%%%%%%%%%%%%%%%%%
\begin{proof}
Given $\lambda, \mu \in \frak{h}^{\ast}_{\op{dom}}$ with
integral difference let
$D^{\lambda}_{\mu} \in D^{\lambda} \op{-Mod_{\DC}-} D^{\mu}$ be the bimodule
$D^{W_{\lambda} \cap W_{\mu}}.$ The preceding considerations show that the
diagram
$$\begin{array}{lcr}
\op{Hom}_{\frak{g}-D} (F,T^{\mu}_{\lambda}K)& \ra
&\op{Hom}_{\frak{g}-D} (T_{\mu}^{\lambda} F, K)\\
\;\;\;\;\;\;\;\;\;\;\;\;\;\;\;\;\;\;\da&&\da\;\;\;\;\;\;\;\;\;\;\;
\;\;\;\;\;\;\;\\
\op{Hom}_{D^{\mu}-D} (\Bbb{V}F, \Bbb{V}
T^{\mu}_{\lambda} K) & 
&\op{Hom}_{D^{\lambda}-D} (\Bbb{V} T^{\lambda}_{\mu}
F, \Bbb{V} K)\\
\;\;\;\;\;\;\;\;\;\;\;\;\;\;\;\;\;\;\da&&\da\;\;\;\;\;\;\;\;\;
\;\;\;\;\;\;\;\;\;\\
 \op{Hom}_{D^{\mu}-D} 
(\Bbb{V}F, \op{Hom}_{D^{\lambda}}(D^{\lambda}_{\mu},
\Bbb{V} K))& \ra
&\op{Hom}_{D^{\lambda}-D} 
(D^{\lambda}_{\mu} \otimes_{D^{\mu}}
\Bbb{V} F, \Bbb{V} K)
\end{array}$$
%%%%%%%%%%%%%%%%%%%%%%%%%%%%%%%%%%%%%%%%%%%%%%%%%%%%%%%%%%%%%%%%%%%%%%
commutes, if we define both lower verticals by
the isomorphisms just introduced and the horizonals by the
adjunctions. In this diagram all morphisms with the possible
exception of both upper verticals are obviously isomorphisms.
Thus if the right upper vertical is an isomorphism, then the left
upper vertical as well.
If in other words our claim holds for
$K$, then it also holds for $T_{\lambda}^{\mu} K$.
Thus it suffices to check it for $K$ a deformed simple Verma.
Working down through a Verma flag, we may even assume
$F$ to be a direct sum of copies of this simple Verma.
In this case the first claim is obvious. The second claim is shown in
the same way.
\end{proof}

\section{Geometrical arguments}

\begin{Remark}\label{AGSU}
  Let $\op{gMod-}\! A$ denote the category of graded right
  modules over a graded ring
$A.$ Let $\op{Der}_G(X)$ resp.\
$\op{Der}_G^+(X)$ denote the equivariant resp.\ bounded below equivariant derived category
corresponding to a complex algebraic variety $X$ 
with the action of a complex algebraic group
$G$ and let $\op{Der}_G(\cal{F},\cal{G})$ denote the morphisms in these
categories, as explained in \cite{BeLu}.
Here and in
 what follows we will always take cohomology with complex coefficients. 
 \end{Remark}
   
 \begin{Remark}
Let
$X$ be a complex algebraic variety with the action of an algebraic group
$B.$
Let $X=\coprod_{a\in A}X_a$ be a stratification into
irreducible locally closed smooth
$B$-stable subvarieties such that the closure of each
stratum is a union of strata.
Let $|a|$ denote the dimension of $X_a$ and
$\cal{C}_a=\underline{X_a}[|a|]$ 
 the ``constant perverse sheaf'' in $\op{Der}_B(X_a).$
Let further $j_a:X_a\hra X$ denote the inclusion.
Let now
$\mathcal{F}, \mathcal{G} \in
\op{Der}_{B} (X)$ be given with the property, that for all
$a\in A$ we have finite direct sum decompositions
\begin{displaymath}
\begin{array}{ccl}
j^{\ast}_{a} \mathcal{F} & \cong & 
\bigoplus_{\nu} f^\nu_{a} \mathcal{C}_a [\nu]\\[2mm]
j^!_{a} \mathcal{G} & \cong & 
\bigoplus_{\nu} g^{\nu}_{a} \mathcal{C}_a [\nu]
\end{array}
\end{displaymath}
in $\op{Der}_{B} (X_a)$ for suitable
$f^{\nu}_a,g^{\nu}_a \in \Bbb{N}.$
If under these assumptions we also have
$f^\nu_{a}= 0=g^\nu_{a}$ for $\nu+|a|$
odd and $H^\nu_B(X_a)= 0$ for $\nu$ odd,
then taking the equivariant hypercohomology 
$\Bbb{H}_B$ induces for all $\ast$  an injection
$$\op{Der}_B(\cal{F},\cal{G}[\ast])\hra 
\op{Hom}(\Bbb{H}_B\cal{F}, \Bbb{H}_B\cal{G})$$
and the dimensions of the homogeneous components on the
left are given by the formula
$$\dim\op{Der}_B(\cal{F},\cal{G}[n])=
\sum_{\nu-\mu+k=n,\; a\in A} f^\nu_a g^\mu_a \dim H^k_B(X_a)$$
The proof is completely analogous to the proof of
Proposition 3 on page 404 of \cite{So-L}
and we shall not repeat it here.
\end{Remark}
%%%%%%%%%%%%%%%%%%%%%%%%%%%%%%%%%%%%%%%%%%%%%%%%%%%%%%%%%%%%%%%%%%%%%%%%%%%%%
\begin{Remark}
Let $G \supset P=P_\iota\supset B \supset T$
be a semisimple complex algebraic group, a parabolic,
a Borel and a maximal torus.
Let $W_\iota\subset W$ be the Weyl groups of $P\subset G$ and
$L\supset T$ the Levi of $P$ above $T.$
We let $B\times P$ act on $G$ by the
rule $(b,p)g=bgp^{-1}.$
From \cite{BeLu} we deduce that the equivariant cohomology
$H^{\ast}_{B\times P}(G)$ becomes under restriction a quotient of
$H^\ast_{B \times P} (\op{pt}) = H^\ast_{T\times L} (\op{pt}) =
  R\otimes_{\Bbb{C}} R^\iota$ for
$R=\DC[\op{Lie}T]$ the ring of regular functions on
$\op{Lie}T$, graded by the condition that linear forms should be
homogeneous of degree two,
and $R^\iota$ the invariants of $W_\iota$ in $R.$
Again using \cite{BeLu} we  get in this way 
a canonical isomorphism
\begin{displaymath}
c : R \otimes_{R^W} R^\iota\overset{\sim}{\rightarrow} 
H^\ast_{B\times P} (G)
\end{displaymath}
Thus
the equivariant hypercohomology
\begin{displaymath}
\Bbb{H}^\ast_{B \times P} : \op{Der}_{B \times P}^+ (G) 
\ra \op{gMod-} H^\ast_{B \times P} (\op{pt})
\end{displaymath}
defines under our identification of the equivariant cohomology ring
and the identification
$\op{Der}_{B\times P}^+ (G)\cong \op{Der}_B^+(G/P)$  a functor
to $\Bbb{Z}$-graded $R$-$R^\iota$-bimodules
\begin{displaymath}
\Bbb{H}_B=\Bbb{H}_B^\ast: \op{Der}^+_{B} (G/P) \ra R\op{-gMod-} R^\iota
\end{displaymath} 
%%%%%%%%%%%%%%%%%%%%%%%%%%%%%%%%%%%%%%%%%%%%%%%%%%%%%%%%%%%%%%%%%%%%%%%%%%%
Now we consider in $\op{Der}_{B}^{+} (G/P)$ for $x \in W/W_\iota$
the intersection cohomology complex
$\cal{IC}_x$ of the closure of ${B x P/P}.$ Let
$\cal{C}_y $ be the constant perverse sheaf
on $B
y B/P,$ which is concentrated in degree $-l(y)$
as a complex of ordinary sheaves, and let $j_y: By B/P
  \hookrightarrow G/P$ denote the embedding.
\end{Remark}
\begin{Satz}\label{VTM}
The functor $\Bbb{H}_B$ is fully faithful on morphisms $\cal{IC}_x \ra
j_{y\ast}\cal{C}_y [n]$ and  $j_{x!}\cal{C}_y  \ra \cal{IC}_y[n]$
and $\cal{IC}_x \ra\cal{IC}_y[n]$
in $\op{Der}_{B}^+ (G/P).$
\end{Satz}
%%%%%%%%%%%%%%%%%%%%%%%%%%%%%%%%%%%%%%%%%%%%%%%%%%%%%%%%%%%%%%%%%%%%%%
\begin{proof}
In \cite{Gi} the corresponding statement is proven for
nonequivariant cohomology and the case $\cal{IC}\ra \cal{IC},$
but in a more general setup.
In \cite{So-L}, Proposition 2, page 402
this is proven for homomorphisms
$\cal{IC}_x \ra\cal{IC}_y[n]$ and $P=B$.
I will now explain in which sense the proof given there
up to some rather minor modifications also proves this more general case.
First we restrict to the case $P=B.$
By Lemma 6 on page 405 of  \cite{So-L}
in connection with \ref{AGSU} the functor of the lemma is faithful
and the dimensions of the $\op{Hom}$ spaces in question are known.
By \cite{SoBi} however we also know the dimensions of the
Hom spaces in the image and thus we may finish the argument
with a comparision of dimensions.
More precisely we get with \ref{AGSU} the formula
\begin{displaymath}
\op{dim}_{\Bbb{C}} 
\op{Der}_{B} (\mathcal{I}\mathcal{C}_x, j_{y\ast} \mathcal{C}_y [n])
=\sum_{k+i = n} n^i_{y,x} \op{dim}_{\Bbb{C}} H^k_{B} (ByB/B)
\end{displaymath}
%%%%%%%%%%%%%%%%%%%%%%%%%%%%%%%%%%%%%%%%%%%%%%%%%%%%%%%%%%%%%%%%%%%%%
Here the $n^i_{y,x}$ are the coefficients of
the Kazhdan-Lusztig polynomials
and we have more precisely
\begin{displaymath}
\sum_{y,i} n^i_{y,x} q^{-i/2} \tilde{T}_{y} = C^{\prime}_{x}
\end{displaymath}
in Lusztig's notation alias
$\sum_{y,i} n^i_{y,x} v^i H_y = \underline{H}_x$ in the notations of
\cite{So-K}.
On the other hand in  \cite{So-L},
Lemma 5, page 402
it is shown for $P=B$, that the $\Bbb{H}_B \cal{IC}_{x}$
are just the special bimodules
  $$\Bbb{H}_B \cal{IC}_{x}\cong B_x$$ 
which I consider in \cite{HCH} and \cite{SoBi}.
Now we recall the graded bimodule
$R_y$ from \cite{SoBi}, which is free of rank one
from the left and the right with the same generator
$1_y$ in degree zero and the property
$r1_y=1_y r^y$ for $r^y=y^{-1}(r),$ and we recall its shifted versions
 $\Delta_y=R_y[-l(y)]$ and
  $\nabla_y=R_y[l(y)].$
  %%%%%%%%%%%%%%%%%%%%%%%%%%%%%%%%%%%%%%%%%%%%%%%%%%%%%%%%%%%%%%%%%%%%%%%%%%%%
One shows easily
$\Bbb{H}_Bj_{y\ast}\cal{C}_y\cong \nabla_y$ and
$\Bbb{H}_Bj_{y!}\cal{C}_y\cong \Delta_y.$
Thus we need to establish the equality of dimenisons
\begin{displaymath}
\op{dim}_{\Bbb{C}} \op{Der}_{B} 
(\mathcal{I}\mathcal{C}_x, j_{y\ast} \mathcal{C}_y[n])
= \op{dim}_{\Bbb{C}} \op{gMod}_{R-R} 
(B_x, \nabla_y
[n]).
\end{displaymath}
But by \cite{SoBi}, Theorem 5.15 the space
$
\op{Mod}_{R-R} (B_x, \nabla_y)
$ 
is graded free as a right $R$-module,
and if we let $h^i_{y,x}$ be the number of generators
needed in degree $i$,
then Theorem 5.3 of   \cite{SoBi} gives in the Hecke algebra
\begin{displaymath}
C^{\prime}_{x} = \sum_{y \in W} h^i_{y,x} q^{-1/2} \tilde{T}_y 
= \sum_{y \in W} h^i_{y,x}
v^i H_y
\end{displaymath}
%%%%%%%%%%%%%%%%%%%%%%%%%%%%%%%%%%%%%%%%%%%%%%%%%%%%%%%%%%%%%
in the notations of Lusztig resp.\ of \cite{So-K}.
In other words we get $h^i_{y,x} = n^i_{y,x},$ and since
$H^{\ast}_{B} (ByB/B) \cong R$
we deduce the claimed equality of dimensions in every degree.
The second case follows dually and thus in case
$P=B$ we have completely established the Lemma.
In general full faithfulness of our functor is deduced in the same way,
but for the equality of dimensions we need to work a little more.
Here we only treat the cases $\mathcal{I} \mathcal{C} \rightarrow
\mathcal{I} \mathcal{C} $ and $\mathcal{I}\mathcal{C} \rightarrow \mathcal{C}$,
the remaining case is dual.
Let $\pi : G/B \twoheadrightarrow G/P$ denote the projection, so
that we get $\Bbb{H}_B \pi_\ast \mathcal{G} \cong \op{res}_{R-R}^{R-R^\iota}
\Bbb{H}_B \mathcal{G}$
and $\mathcal{H}_B\pi^\ast \mathcal{F} \cong \Bbb{H}_B\mathcal{F}
\otimes_{R^\iota} R.$
This leads to a commutative diagramm
\begin{displaymath}
\begin{array}{ccc}
\op{Der}_B (\mathcal{F}, \pi_\ast \mathcal{G} \left[ n\right]) & 
\overset{\sim}{\longrightarrow} & \op{Der}_B (\pi^\ast \mathcal{F}, 
\mathcal{G} \left[ n 
\right])\\
\downarrow & & \downarrow \\
\op{gMod}_{R-R^\iota} (\Bbb{H}_B \mathcal{F}, \Bbb{H}_B \pi_\ast 
\mathcal{G} \left[ n
\right]) & & \op{gMod}_{R-R}
(\Bbb{H}_B \pi^\ast \mathcal{F}, \Bbb{H}_B \mathcal{G} \left[ n
\right])\\
\parallel & &\parallel\\
\op{gMod}_{R-R^\iota} (\Bbb{H}_B \mathcal{F}, \op{res}_{R-R}^{R-R^\iota} 
\Bbb{H}_B \mathcal{G}
\left[ n \right]) &
\overset{\sim}{\longrightarrow} & \op{gMod}_{R-R} (\Bbb{H}_B \mathcal{F} 
\otimes_{R^\iota} R,
\Bbb{H}_B \mathcal{G} \left[ n \right])
\end{array}
\end{displaymath}
%%%%%%%%%%%%%%%%%%%%%%%%%%%%%%%%%%%%%%%%%%%%%%%%%%%%%%%%%%%%%
and with the right upper vertical the left upper vertical must be
an isomorphism, too.
Thus the cases $\mathcal{I}\mathcal{C}
\rightarrow \mathcal{I}\mathcal{C}$
and $\mathcal{I}\mathcal{C} \rightarrow \mathcal{C}$
follow for general $P$ from the case $P=B$.
\end{proof}
%%%%%%%%%%%%%%%%%%%%%%%%%%%%%%%%%%%%%%%%%%%%%%%%%%%%%%%

\section{Singular bimodules}
\begin{Remark}
Let $\mathcal{W}$ be a finite group of automorphisms of a
finite dimensional affine space
$E$ over $\Bbb{Q}$, which is generated by reflections,
and let $\mathcal{S} \subset \mathcal{W}$
be a choice of simple reflections.
Let $R$ denote the regular functions on the space of translations, graded
by the rule, that linear functions are homogeneous of degree two.
Then by \cite{SoBi} there exist well-defined up to
isomorphism $\Bbb{Z}$-graded
$R$-bimodules $B_x = B_x (\mathcal{W}) = B_x (\mathcal{W},\mathcal{S}, E)
\in R\op{-gMod-}R$ such that we have
\begin{enumerate}
\item The $B_x$ are indecomposable.
\item
For $e$ the neutral element we have $B_e=R.$
\item If $s \in \mathcal{S} $ is a simple reflection with
$x s> x$, then there is a decomposition
\begin{displaymath}
B_x \otimes_{R^s} R [1] \cong B_{xs} \oplus \bigoplus_{l(y) \leq l(x)}
{m(y)}B_y
\end{displaymath}
for suitable multiplicities $m(y) \in \Bbb{N}$.
\end{enumerate}
Following \cite{SoBi}
the rings of endomorphisms of degree zero of
these bimodules consist just of scalars,
in particular our bimodules stay indecomposable when we extend scalars.
Now let $\mathcal{S}_{\iota} \subset \mathcal{S}$ be a subset of the
set of simple reflections,
$\mathcal{W}_\iota = \langle \mathcal{S}_\iota
\rangle \subset \mathcal{W}$
the subgroup generated by it, $w_\iota\in{\mathcal{W}}_\iota$ its longest
element and $R^\iota$  the
subring of ${\mathcal{W}_\iota}$-invariants. Then under the same
assumptions we claim:
\end{Remark}
%%%%%%%%%%%%%%%%%%%%%%%%%%%%%%%%%%%%%%%%%%%%%%%%%%%%%%%%%%%%%%%%%%%%%%%%%%%%%%%%%%%%%%
\begin{Lemma}
For every coset $\bar{x} \in \mathcal{W}/\mathcal{W}_\iota$ there exists
one and only one $\Bbb{Z}$-graded $R$-$R^\iota$-bimodule
$B_{\bar{x}}^\iota =
B^\iota_{\bar{x}} (\mathcal{W},\mathcal{W}_\iota) \in R
\op{-gMod-}R^\iota$ with the property that for $x \in \mathcal{W}$ the
longest representative of the coset $\bar{x}$ we have
\begin{displaymath}
\op{res}^{R-R^\iota}_{R-R} B_x \cong \bigoplus_{z \in 
\mathcal{W}_\iota} B^\iota_{\bar{x}}
[l(w_\iota)- 2l (z)]
\end{displaymath}
\end{Lemma}
%%%%%%%%%%%%%%%%%%%%%%%%%%%%%%%%%%%%%%%%%%%%%%%%%%%%%%%%%%%%%%%%%%%%%%%%%%%%%%%%%%%%%
\begin{proof}
Without restriction of generality we may assume that
$\mathcal{W}$ admits only one fixed point.
Since by assumption it is a rational and thus crystallographic
reflection group, we then
find $G \supset B \supset T$ a complex semisimple algebraic group
$G$ with Borel $B$ and maximal
torus $T$ and Coxeter system $(\mathcal{W},
\mathcal{S})$.
If we identify in a $\mathcal{W}$-equivariant way
$H^2_{T} (\op{pt};\Bbb{Q})$ and the homogeneous component
$R^2$ of $R$, then as we discussed already 
in the proof of \ref{VTM}  there exists an isomorphism of
$\Bbb{Z}$-graded $R$-bimodules
\begin{displaymath}
B_x \cong \Bbb{H}_B \mathcal{I}\mathcal{C} (\overline{BxB / B})
\end{displaymath}
%%%%%%%%%%%%%%%%%%%%%%%%%%%%%%%%%%%%%%%%%%%%%%%%%%%%%%%%%%%%%%%%%%%%%%%%%%%
If $P=P_\iota$
is a parabolic with $G \supset P \supset B$, then the decomposition theorem
of \cite{BeLu} applied to the projection $p: G/B \twoheadrightarrow G/P,$
shows for $x$ maximal in its $\mathcal{W}_{\iota}$-coset the
existence of a decomposition
\begin{displaymath}
p_\ast \mathcal{I}\mathcal{C} (\overline{BxB/B}) \cong 
\bigoplus_{z\in \mathcal{W}_\iota} \mathcal{I}\mathcal{C}
(\overline{Bx P/P}) [l(w_\iota)-2l(z)]
\end{displaymath}
in $\op{Der}_B^+ (G/P)$. 
%%%%%%%%%%%%%%%%%%%%%%%%%%%%%%%%%%%%%%%%%%%%%%%%%%%%%%%%%%%%%%%%%%%%%%%%%%%%%
On the other hand we have
$\Bbb{H}_B \circ p_\ast=  \op{res}^{R-R^\iota}_{R-R}   \circ\Bbb{H}_B$
and the $\Bbb{H}_B\mathcal{I}\mathcal{C}
(\overline{Bx P/P})$ are indecomposable as graded $R$-$R^\iota$-bimodules,
since by \ref{VTM} the scalars are their only endomorphisms of degree
$\leq 0$.
\end{proof}
%%%%%%%%%%%%%%%%%%%%%%%%%%%%%%%%%%%%%%%%%%%%%%%%%%%%%%%%%%%%%%%%%%%%%%%%%%%%%%
%%%%%%%%%%%%%%%%%%%%%%%%%%%%%%%%%%%%%%%%%%%%%%%%%%%%%%%%%%%%%%%%%%%%%%%%%%%%%%%
\section{The bimodules for tilting modules}
\begin{Remark}
Given $y\in W$ let $ \hat{S}_{{y}}$ denote the bimodule, which
from the left is free over $\hat{S}$ of rank one with basis 
say $1_y,$ but from the right has the action
$r 1_y = 1_y r^y$ of $\hat{S}.$
Given a bimodule $B$ for two commutative rings let
$\tilde{B}$ denote the bimodule which we get by interchanging
the right and the left action.
\end{Remark}
%%%%%%%%%%%%%%%%%%%%%%%%%%%%%%%%%%%%%%%%%%%%%%%%%%%%%%%%%%%%%%%%%%%%%
\begin{Satz}\label{BiKi}
Let $\lambda \in \frak{h}^*_{\op{dom}}$ and let $W_\lambda
\subset W_{\bar{\lambda}} \subset W$ be as in \ref{WLL}.
Given $\bar{x} \in W_{\bar{\lambda}} / W_\lambda$ we have 
for the deformation of the indecomposable tilting module with highest
weight $w_{\bar{\lambda}} \bar{x} \cdot \lambda$ the formula
\begin{displaymath}
\Bbb{V} K_{\hat{S}} (w_{\bar{\lambda}} \bar{x} \cdot \lambda)
\cong \widetilde{B}^\lambda_{\bar{x}} 
\otimes_S \hat{S}_{w_{\bar{\lambda}}}
\end{displaymath}
\end{Satz}
\begin{Remarkc}
Our bimodules $B^\iota_{\bar{x}}$ are graded free of finite rank over
$R$, thus $\hat{R} \otimes_{R} B^\iota_{\bar{x}}$
is just the completion along the grading of
our original bimodule $B^\iota_{\bar{x}}$ and in particular 
this completion admits a right action of $\hat{R}^\iota$.
\end{Remarkc}
\begin{proof}
As is well-known \korr{and explained in Remark 7.2.2 of \cite{So-K},}
an indecomposable tilting module stays
indecomposable upon translation out of the walls. More precisely
for $\lambda, \mu \in \frak{h}^*_{\op{dom}}$ with
$\lambda + X = \mu +X$ and $W_\mu =1$ and $x \in W_{\bar{\lambda}}$
maximal in its coset $xW_{\lambda}$ we have
\begin{displaymath}
T^\mu_\lambda K (w_{\bar{\lambda}} x \cdot \lambda) =
K (w_{\bar{\mu}} x \cdot \mu)
\end{displaymath}
Since $T^\lambda_\mu T^\mu_\lambda$ is a sum of
$|W_\lambda|$ copies of the identity functor by \cite{BG}, we may for the
proof restrict to the case
$\lambda$ regular.
If $x =st \ldots r$ is a reduced decomposition by simple
reflections of
$W_{\bar{\lambda}}$, we may characterize
$K_{\hat{S}} (w_{\bar{\lambda}}
x \cdot  \lambda)$
inductively as the indecomposable summand of
\begin{displaymath}
\vartheta_r \ldots \vartheta_t \vartheta_s 
\Delta_{\hat{S}} (w_{\bar{\lambda}}
\cdot \lambda)
\end{displaymath}
not isomorphic to any
$K_{\hat{S}} (w_{\bar{\lambda}}y \cdot \lambda)$
for $y < x$.
Applying $\Bbb{V}$ we get from this the
indecomposable summand of
\begin{displaymath}
\hat{S} \otimes_{\hat{S}^r} 
\hat{S} \ldots \otimes_{\hat{S}^t} \hat{S}
\otimes_{\hat{S}^s} \hat{S}_{w_{\bar{\lambda}}}
\end{displaymath}
which didn't appear already before.
But by definition of our special bimodules this is precisely
$\widetilde{B}_x \otimes_S \hat{S}_{w_{\bar{\lambda}}}$.
\end{proof}
%%%%%%%%%%%%%%%%%%%%%%%%%%%%%%%%%%%%%%%%%%%%%%%%%%%%%%%%%%%%%%%%%%%%%%%%%%
\section{Restricting the group action}
%%%%%%%%%%%%%%%%%%%%%%%%%%%%%%%%%%%%%%%%%%%%%%%%%%%%%%%%%%%%%%%%%%%%%%%%%
\begin{Remark}
For $G$ a complex connected algebraic group
we let $A_G=H^\ast(BG)$ be the cohomology ring of its 
classifying space.
If $X$ is a complex algebraic
$G$-variety and  $\mathcal{F}, \mathcal{G} \in
\op{Der}^+_{G} (X)$ are objects of the equivariant derived category,
we may form the graded $A_G$-module
$$\op{Der}_G (\mathcal{F}, \mathcal{G} [\ast])
=\bigoplus_n \op{Der}_G (\mathcal{F}, \mathcal{G} [n])$$
\end{Remark}
\begin{Proposition}\label{GruWe}
Let $G \supset H$ be a connected complex algebraic group
and a connected closed subgroup.
Let $X$ be an algebraic $G$-variety and let
$\mathcal{F}, \mathcal{G} \in
\op{Der}_{G} (X)$ be constructible complexes.
If $\op{Der}_G (\mathcal{F}, \mathcal{G} [\ast])$ is graded free
over $A_G,$
then the obvious map induces a bijection
\begin{displaymath}
A_H \otimes_{A_G} 
 \op{Der}_G (\mathcal{F},\mathcal{G} 
[\ast]) \overset{\sim}{\rightarrow} 
\op{Der}_H (\mathcal{F},\mathcal{G} [\ast])
\end{displaymath}
\end{Proposition}
\begin{proof}
We consider the constant map
$k: X \ra \op{pt}$ and the fully faithful functor
$
\gamma_G : \op{Der}_G^{\op{c}} (\op{pt}) \ra A_G\op{-dgDer} 
$
from \cite{BeLu}, 12.4.6 and recall
$$
\op{Der}_G (\mathcal{F}, \mathcal{G} [\ast]) = H^\ast \gamma_G k_{\ast}
\op{Hom}(\mathcal{F},\mathcal{G})
$$
where we form $\op{Hom}(\mathcal{F}, \mathcal{G})$
in $\op{Der}_G^+ (X)$ and
$k_{\ast}$ means the direct image landing in $\op{Der}_G^+ (\op{pt})$.
If this now is a free $A_G$-module,
then $\gamma_G k_\ast \op{Hom} (\mathcal{F},\mathcal{G})$
is already quasiisomorphic to its cohomology and this cohomology is
homotopy projective in
$A_G\op{-dgMod}.$
With the derived functor
$A_H \otimes^{L}_{A_G}: A_G\op{-dgDer}\ra A_H\op{-dgDer}$
we have by \cite{BeLu}, 12.7.1 furthermore canonically
\begin{displaymath}
\left(A_H \otimes^{L}_{A_G}\right) \circ\gamma_G
= \gamma_{H} \circ \op{res}^H_G
\end{displaymath}
and for homotopy projective objects
$M \in A_G\op{-dgMod}$ we have in addition
\begin{displaymath}
A_H \otimes^L_{A_G} M
= A_H \otimes_{A_G} M
\end{displaymath}
Since $k_{\ast}$ and $\op{Hom}$ commute with the restriction
of the group action, this shows the Proposition.
\end{proof}
%%%%%%%%%%%%%%%%%%%%%%%%%%%%%%%%%%%%%%%%%%%%%%%%%%%%%%%%%%%%%%%%%%%%%%%%%%%%%%%%%%%%%

\section[Proof of the main theorem 4.4]{Proof of the main theorem \protect\ref{MT}}
%%%%%%%%%%%%%%%%%%%%%%%%%%%%%%%%%%%%%%%%%%%%%%%%%%%%%%%%%%%%%%%%%%%%%%%%%%%%%%%%%%
\begin{Remark}
Given a ring $R,$  a ring homomorphism
$R \ra \Bbb{C} [[v]],$
three $R$-modules $H,H',H''$ and an $R$-bilinear map
\begin{displaymath}
\varphi : H \times H^\prime \rightarrow H^{''}
\end{displaymath}
we may introduce on $\Bbb{C} [[v]]\otimes_R H$ a filtration by the rule
\begin{displaymath}
F^i (\Bbb{C} [[v]] \otimes_R H) = \left\{
h \left| \begin{array}{l}
\varphi (h,h^\prime)\in v^i \Bbb{C} [[v]] \otimes_R H^{''}\\
\text{for all } h^\prime \in \Bbb{C} [[v]] \otimes_R H^\prime
\end{array} \right. \right\}
\end{displaymath}
and get also an induced filtration on
$\Bbb{C} \otimes_R H,$
whose subquotients we denote $\bar{F}^i(\Bbb{C} \otimes_R H)$ .
\end{Remark}
%%%%%%%%%%%%%%%%%%%%%%%%%%%%%%%%%%%%%%%%%%%%%%%%%%%%%%%%%%%%%%%%%%%%%%%%%%
\begin{Remark}\label{NZ}
If for example $\hat{S}$ is the completion of
$S = \mathcal{O} (\frak{h}^*)$
along the natural grading and $\hat{S}
\twoheadrightarrow \Bbb{C} [[v]]$ is the restriction to the line
$\DC\rho$ as in \ref{DAF},
then the pairing given by composition
\begin{displaymath}
\begin{array}{r}
\op{Hom}_{\frak{g} - \hat{S}} 
(\Delta_{\hat{S}} (\lambda), K_{\hat{S}} (\mu)) \times 
\op{Hom}_{\frak{g}-\hat{S}} (K_{\hat{S}} (\mu), 
\nabla_{\hat{S}} (\lambda))\\[2mm] \rightarrow
\op{Hom}_{\frak{g}-\hat{S}} (\Delta_{\hat{S}} (\lambda),
\nabla_{\hat{S}} (\lambda))
\end{array}
\end{displaymath}
leads to the Andersen filtration \ref{DAF} on the spaces
$\op{Hom}_{\frak{g}} (\Delta (\lambda), K (\mu)),$
which by \ref{BCH} may be identified with
$
\op{Hom}_{\frak{g}-\hat{S}} (\Delta_{\hat{S}} (\lambda), 
K_{\hat{S}} (\mu)) \otimes_{\hat{S}}
\Bbb{C}$ .
Certainly there also exists
$p\in \hat{S}$ such that $\Bbb{V}$ induces an isomorphism
\begin{displaymath}
\op{Hom}_{\frak{g}-\hat{S}} (\Delta_{\hat{S}} (\lambda), 
\nabla_{\hat{S}} (\lambda))
\overset{\sim}{\rightarrow} p \op{Hom}_{\hat{S}} 
(\Bbb{V} \Delta_{\hat{S} }(\lambda),
\Bbb{V} \nabla_{\hat{S}} (\lambda))
\end{displaymath}
and with this $p$ we may reformulate our pairing as in
the pairing given by the composition
\begin{displaymath}
\begin{array}{r}
\op{Hom}_{\hat{S}^{\lambda} -\hat{S}} 
(\Bbb{V} \Delta_{\hat{S}} (\lambda),
\Bbb{V} K_{\hat{S}}(\mu)) \times  
\op{Hom}_{\hat{S}^{\lambda} - \hat{S}}
(\Bbb{V} K_{\hat{S}} (\mu), \Bbb{V} \nabla_{\hat{S}} (\lambda))\\[2mm]
\rightarrow p \op{Hom}_{\hat{S}^{\lambda}-\hat{S}} 
(\Bbb{V} \Delta_{\hat{S}}
(\lambda), \Bbb{V} \nabla_{\hat{S}} (\lambda))
\end{array}
\end{displaymath}
Here a possible $p$ may be determined by the condition that
our pairing in case
$\lambda = \mu$ must lead to a surjection.
%%%%%%%%%%%%%%%%%%%%%%%%%%%%%%%%%%%%%%%%%%%%%%%%%%%%%%%%%%%%%%%%%%%
Now we change parameters, choose $\lambda \in \frak{h}^*_{\op{dom}}$ and
put $\lambda_{\bar{x}} = w_{\bar{\lambda}} \bar{x} \cdot \lambda$
for $\bar{x} \in W_{\bar{\lambda}}/ W_{\lambda}$.
To simplify we further introduce a variant
$\widetilde{\Bbb{V}}$ of $\Bbb{V}$ by putting
\begin{displaymath}
\widetilde{\Bbb{V}} M = \hat{S}_{w_{\bar{\lambda}}} 
\otimes_{\hat{S}}
\widetilde{\Bbb{V}M}
\end{displaymath}
such that \ref{BiKi} becomes
$
\widetilde{\Bbb{V}} K_{\hat{S}} (\lambda_{\bar{x}}) 
\cong \hat{B}^\lambda_{\bar{x}}
,$
the hat meaning completion along the grading.
With less effort one may also check $$\widetilde{\Bbb{V}}
\Delta_{\hat{S}} (\lambda_{\bar{y}})
\cong \widetilde{\Bbb{V}} \nabla_{\hat{S}} 
(\lambda_{\bar{y}}) \cong \hat{S}^\lambda_{\bar{y}}$$
in $\hat{S}\op{-mod-}\hat{S}^\lambda$, 
where again we mean the bimodule which is
$\hat{S}$ from the left but has the
$\bar{y}$-twisted action
$r 1_y = 1_y r^y$ of $\hat{S}^\lambda$ from the right.
If we replace $\Bbb{V}$ by $\widetilde{\Bbb{V}}$,
we thus obtain up to a twist of the right
$\hat{S}$-action by $w_{\bar{\lambda}}$ the pairing 
\begin{displaymath}
  \begin{array}{r}
\op{Hom}_{\hat{S}-\hat{S}^\lambda} (\hat{S}^\lambda_{\bar{y}},
\hat{B}^\lambda_{\bar{x}})\times \op{Hom}_{\hat{S}-\hat{S}^\lambda}
(\hat{B}^\lambda_{\bar{x}}, \hat{S}^\lambda_{\bar{y}})
\\[2mm] 
\rightarrow p_1 
\op{Hom}_{\hat{S}-\hat{S}^\lambda}({\hat{S}^\lambda_{\bar{y}}}, 
\hat{S}^\lambda_{\bar{y}})
\end{array}
\end{displaymath}
of $\hat{S}$-modules and our filtration corresponds to the
filtration we get here when we change
$\hat{S} \twoheadrightarrow \Bbb{C} [[v]]$ by twisting
it with $w_{\bar{\lambda}}$.
Here $p_1$ denotes the image of
$p$ under  $w_{\bar{\lambda}}.$
%%%%%%%%%%%%%%%%%%%%%%%%%%%%%%%%%%%%%%%%%%%%%%%%%%%%%%%%%%%%%%%%%%%%%%%%%%%%%
Since the choice of $p_1$ is only sensible up to units of
$\hat{S}$,  we may choose  $p_1$ already before completion and the
corresponding pairing
 ``before completion''
\begin{displaymath}
  \begin{array}{r}
\op{Hom}_{S-S^\lambda} (S^\lambda_{\bar{y}}, 
B^\lambda_{\bar{x}}) \times \op{Hom}_{S-S^\lambda}
(B^\lambda_{\bar{x}}, S^\lambda_{\bar{y}}) \\[2mm] \rightarrow
p_1 \op{Hom}_{S-S^\lambda} (S^\lambda_{\bar{y}}, S^\lambda_{\bar{y}})
\end{array}
\end{displaymath}
leads to the same filtered $\Bbb{C}$-space in the end.
This pairing we now interpret geometrically.
\end{Remark}
  
\begin{Remark}
If $X \subset \frak{h}^*$ denotes the lattice of integral weights,
we find a pair
$G^\vee \supset T^\vee$ consisting of a reductive
connected complex algebraic group with a maximal torus
such that $X = X(T^\vee)$
is its group of one-parameter subgroups and that for the
Weyl group we have $W( G^\vee, T^\vee)=W_{\bar{\lambda}}.$
In $G^\vee$ we then choose a Borel $B^\vee$ for $\mathcal{S}
\cap W_{\bar{\lambda}}$ and a parabolic
$P^\vee \supset B^\vee$ for $W_\lambda$.
If now $$\mathcal{I} \mathcal{C}_{\bar{x}}
= \mathcal{I}\mathcal{C}
(\overline{B^\vee \bar{x} P^\vee/P^\vee})$$
denotes the intersection homology complex of the corresponding Schubert variety and
$\mathcal{C}_{\bar{y}}$ the constant perverse sheaf on
$B^\vee\bar{y} P^\vee/P^\vee$,
we have $B^\lambda_{\bar{x}} \cong
\Bbb{H}_{B^\vee} \mathcal{I}\mathcal{C}_{\bar{x}}$
and our pairing ``before completion'' from the end of the previous
remark can be interpreted with the help of \ref{VTM} as
the pairing given by composition
\begin{displaymath}
  \begin{array}{r}
\op{Der} _{B^\vee} (j_{\bar{y}!}\mathcal{C}_{\bar{y}}, 
\mathcal{I}\mathcal{C}_{\bar{x}}[*])
\times \op{Der} _{B^\vee} (\mathcal{I}\mathcal{C}_{\bar{x}}, 
j_{\bar{y}\ast} \mathcal{C}_{\bar{y}}[*])\\[2mm]
\rightarrow \op{Der} _{B^\vee} (j_{\bar{y}!} 
\mathcal{C}_{\bar{y}}, j_{\bar{y}\ast}
\mathcal{C}_{\bar{y}}[*])
\end{array}
\end{displaymath}
%%%%%%%%%%%%%%%%%%%%%%%%%%%%%%%%%%%%%%%%%%%%%%%%%%%%%%%%%%%%%%%%%%%%%%%%%%%%%%%
Here we do not need a $p$-factor on
the right hand side,
since restriction to the big cell shows that for
$\bar{x}=\bar{y}$ our pairing gives a surjection.
The question is thus, which filtered vector space this
pairing of
$A_{B^\vee}$-modules leads to under the homomorphism
$A_{B^\vee} \twoheadrightarrow \Bbb{C} [v]$
coming from the embedding $\Bbb{C}^\times \hookrightarrow T^\vee$
with parameter $w_{\bar{\lambda}}
\rho$.
But by \ref{GruWe} this specialization leads us to the composition pairing
%$\Bbb{C} [v]$
\begin{displaymath}
  \begin{array}{r}
\op{Der} _{\Bbb{C}^\times} (j_{\bar{y}!}\mathcal{C}_{\bar{y}}, 
\mathcal{I}
\mathcal{C}_{\bar{x}}[*]) \times \op{Der} _{\Bbb{C}^\times}
(\mathcal{I}\mathcal{C}_{\bar{x}}, j_{\bar{y}\ast} 
\mathcal{C}_{\bar{y}}[*])\\[2mm]
\rightarrow \op{Der} _{\Bbb{C}^\times} (j_{\bar{y}!}
\mathcal{C}_{\bar{y}},
j_{\bar{y}\ast} \mathcal{C}_{\bar{y}}[*])
\end{array}
\end{displaymath}
Let now $\bar{y}$ denote the point
$\bar{y}P^\vee$ of $G^\vee/P^\vee.$
For a suitable product $U$ of root subgroups of $G^\vee$ the multiplication
$u \mapsto u{\bar{y}}$ defines an embedding
$U \hookrightarrow G^\vee / P^\vee,$ 
whose image is a cell transversal to $B^\vee{\bar{y}} P^\vee / P^\vee$
and is contracted by
$\Bbb{C}^\times$ to  $\bar{y}$.
If we put $Z = U{\bar{y} }\cap \overline{B^\vee\bar{x} P^\vee/P^\vee},$
then $Z$ is contracted by $\Bbb{C}^\times$ to $\bar{y}$, and
if $a : Z \hookrightarrow G^\vee/P^\vee$ denotes the embedding,
the restriction to $Z$ will not change our pairing.
If we now put  $d =\dim B^\vee{\bar{y} }P^\vee/P^\vee$ and let
$i:\op{pt}\hra Z$ be the embedding of $\bar{y}$
and $\underline{\op{pt}}$ the constant sheaf on a point, we get
$
a^* j_{\bar{y}!} \mathcal{C}_{\bar{y}} \cong i_* \underline{\op{pt}}
[d] \cong a^*j_{\bar{y}_*} \mathcal{C}_{\bar{y}}
$ 
and $a^* \mathcal{I}\mathcal{C}_{\bar{x}} \cong \mathcal{I} \mathcal{C}
 [d]$ will be the shifted intersection cohomology complex 
 $\mathcal{I} \mathcal{C}
=\mathcal{I} \mathcal{C}(Z)$ of $Z$
and our pairing gets transformed to the composition pairing
\begin{displaymath}
\begin{array}{r}
\op{Der}_{\Bbb{C}^\times} (i_* \underline{\op{pt}}, 
\mathcal{I}\mathcal{C}[*]) \times \op{Der}_{\Bbb{C}^\times}
(\mathcal{I}\mathcal{C} , i_* \underline{\op{pt}}[*])\\[2mm]
 \rightarrow  \op{Der}_{\Bbb{C}^\times} (i_* \underline{\op{pt}}, 
i_* \underline{\op{pt}}[*])
\end{array}
\end{displaymath}
Now we may identify the first of these paired modules with
$H_{\DC^\times}(i^! \mathcal{I}
\mathcal{C}) $ and the second with the dual of 
$H_{\DC^\times}(i^* \mathcal{I}\mathcal{C}) $ 
and thus our pairing leads to the same filtered space
as the embedding of free
$\Bbb{C}[v]$-modules
$H_{\DC^\times}(i^! \mathcal{I}\mathcal{C})  
\hookrightarrow H_{\DC^\times}(i^* \mathcal{I}
\mathcal{C}) .$
But in this situation the  ``fundamental example''
of section 14 of \cite{BeLu}  just means that the
cokernel of this embedding may be identified with the intersection
cohomology of the projective variety
$\overline{Z} = (Z\setminus \{\overline{y}\})
/\Bbb{C}^\times$ shifted by one,
with $IC (\overline{Z}) [1]$
viewed as a $\Bbb{C} [v]$-module
in such a way, that
$v$ acts as  Lefschetz operator,
thus leading to a short exact sequence of $\DZ$-graded $\DC[v]$-modules 
$$H_{\DC^\times}(i^! \mathcal{I}\mathcal{C})  
\hookrightarrow H_{\DC^\times}(i^* \mathcal{I}
\mathcal{C})\sra IC (\overline{Z}) [1]$$
More precisely, this goes as follows: One starts with the decomposition
of $Z$ into the center of the contraction $\op{pt}$ and its
open complement $Z_0$ and denotes the inclusions by $i$ and $j$ and 
writes the Gysin triangle in the equivariant derived category 
pushed down by a map $p$ to a point
$$p_\ast i_! i^! \mathcal{I}\mathcal{C} \ra p_\ast\mathcal{I}\mathcal{C} \ra 
p_\ast j_* j^* \mathcal{I}\mathcal{C}\stackrel{[1]}{\ra}$$
This 
triangle in $\op{Der}_{\DC^\times}(\op{pt})$ 
can by Theorem 14.2 in loc.cit. 
be identified with a triangle in the derived 
category of dg-modules  $\DC[v]\op{-dgDer}$ over 
$\DC[v]=H^\ast_{\DC^\times}(\op{pt})$
written
$$(\tau_{\geq 0}M)[-1]\ra \tau_{< 0}M\ra M\stackrel{[1]}{\ra}$$
in the notation of loc.cit. However at the end of the proof of
this Theorem following loc.cit. 
Lemma 14.15 this is further rewritten as a short 
exact sequence of graded $\DC[v]$-modules and
$M$ is identified with $IC (\overline{Z}) [1]$ by a remark preceding 
loc.cit. 14.5 labeled 13.4, since in our 
case $G$ is trivial, whereas $(\tau_{\geq 0}M)[-1]$
is identified with the costalk by loc.cit. Theorem 14.2(i) and this
costalk with its cohomology by loc.cit. 14.3(i).
The middle part of our Gysin sequence finally may also be
interpreted as the  stalk at the center of
the contraction, 
$p_\ast \cal{I}\cal{C}\cong i^\ast \cal{I}\cal{C},$
as explained in \cite{Spp}, section 3, and again this stalk 
when written as a dg-module
can be
identified  with its cohomology by \cite{BeLu} 14.3(i').
In this way we see that \cite{BeLu} indeed leads to the short 
exact sequence of $\DC[v]$-modules I claimed.
The hard Lefschetz
from \cite{BBD} 6.2.10 now tells us,
that as a graded module over $\DC[v]$ our
intersection cohomology $IC (\overline{Z})$ is a direct 
sum of truncated polynomial rings graded in such a way 
they are selfdual, in formulas a direct sum of graded 
modules of the form $(\DC[v]/(v^{i+1}))[i].$
Furthermore we know from the description of our sequence by truncation,
say, that its middle module is freely generated in negative degrees and 
its first module is freely generated in positive degrees.
This means that on the level of graded $\DC[v]$-modules, 
our sequence breaks up into a direct sum of copies of sequences 
of the type
$$(\DC[v])[-i]\stackrel{v^i\cdot}{\hra}(\DC[v])[i]\sra (\DC[v]/(v^{i}))[i]$$
with $i> 0.$
From there we see easily, 
that the said filtration on 
$\Bbb{C} \otimes_{\Bbb{C}[v]} 
H_{\DC^\times}(i^! \mathcal{I}\mathcal{C})=
H(i^! \mathcal{I}\mathcal{C}) $
coincides with the filtration given by its
$\Bbb{Z}$-grading, and this by \cite{KL-S} is known to be given by
 Kazhdan-Lusztig polynomials.
More precisely we get
\begin{displaymath}
\begin{array}{ccl}
\bar{F}^i (\Bbb{C} \otimes_{\Bbb{C}[v]} i^! \mathcal{I}\mathcal{C}) &
\cong & H^{-i} (i^! \mathcal{I}\mathcal{C})\\
& \cong & \op{Der} (\mathcal{C}_{\bar{y}} [i],i^!_{\bar{y}} 
\mathcal{I}
\mathcal{C}_{\bar{x}})\\
&\cong & \op{Der} (i_{\bar{y}!} \mathcal{C}_{\bar{y}} [i], 
\mathcal{I}
\mathcal{C}_{\bar{x}})
\end{array}
\end{displaymath}
and this space has the dimension $h^i_{y,x}$ for $y,x$
the longest representatives of
$\bar{y}, \bar{x}$.
Thus this is the dimension of the $i$-th subquotient of the
Andersen filtration on
$
\op{Hom}_{\frak{g}}(\Delta (\lambda_{\bar{y}}), 
K(\lambda_{\bar{x}}))
.$
\end{Remark}

%\bibliographystyle{amsalpha}\bibliography{pub}
%\end{document}

\providecommand{\bysame}{\leavevmode\hbox to3em{\hrulefill}\thinspace}
\providecommand{\MR}{\relax\ifhmode\unskip\space\fi MR }
\providecommand{\MRhref}[2]{%
  \href{http://www.ams.org/mathscinet-getitem?mr=#1}{#2}
}
\providecommand{\href}[2]{#2}

\end{document}